# Local times of multifractional Brownian sheets

MARK MEERSCHAERT[1,*], DONGSHENG WU[2] and YIMIN XIAO[1,**]

[1]*Department of Statistics and Probability, Michigan State University, East Lansing, MI 48824, USA. E-mail:* [*]*mcubed@stt.msu.edu;* [**]*xiao@stt.msu.edu*

[2]*Department of Mathematical Sciences, University of Alabama in Huntsville, Huntsville, AL 35899, USA. E-mail:* *dongsheng.wu@uah.edu*

Denote by $H(t) = (H_1(t), \ldots, H_N(t))$ a function in $t \in \mathbb{R}_+^N$ with values in $(0,1)^N$. Let $\{B^{H(t)}(t)\} = \{B^{H(t)}(t), t \in \mathbb{R}_+^N\}$ be an $(N,d)$-multifractional Brownian sheet (mfBs) with Hurst functional $H(t)$. Under some regularity conditions on the function $H(t)$, we prove the existence, joint continuity and the Hölder regularity of the local times of $\{B^{H(t)}(t)\}$. We also determine the Hausdorff dimensions of the level sets of $\{B^{H(t)}(t)\}$. Our results extend the corresponding results for fractional Brownian sheets and multifractional Brownian motion to multifractional Brownian sheets.

## 1. Introduction

A one-dimensional fractional Brownian motion (fBm) $\xi^\alpha = \{\xi^\alpha(t), t \in \mathbb{R}\}$ with Hurst index $\alpha \in (0,1)$ is a real-valued, centered Gaussian process with covariance function given by

$$\mathbb{E}[\xi^\alpha(s)\xi^\alpha(t)] = \tfrac{1}{2}[|s|^{2\alpha} + |t|^{2\alpha} - |t-s|^{2\alpha}] \qquad \forall s,t \in \mathbb{R}. \tag{1.1}$$

It was introduced, as a moving average Gaussian process, by Mandelbrot and Van Ness [31].

Fractional Brownian motion has interesting properties such as self-similarity of order $\alpha \in (0,1)$, stationary increments and long-range dependence (when $\alpha > 1/2$), which make it a good candidate for modeling different phenomena in, for example, finance and telecommunication. However, this model may be restrictive, due to the fact that all of its regularity and fractal properties are governed by the single Hurst parameter $\alpha$. To model phenomena whose regularity evolves in time, for example, Internet traffic or image processing, Lévy-Véhel and Peltier [28] and Benassi, Jaffard and Roux [7] independently introduced multifractional Brownian motion (mfBm) in terms of moving average representation and harmonizable representation, respectively. Multifractional Brownian motion is governed by a Hurst functional $\alpha(t)$ with certain regularity in place of the constant Hurst parameter $\alpha$ in fBm.

Several authors have studied sample path and statistical properties of mfBm. For example, Benassi, Jaffard and Roux [7] considered the sample path Hölder regularity of mfBm and determined the Hausdorff dimension of its graph. Ayache, Cohen and Lévy-Véhel [2] and Herbin [22] studied the covariance structure of mfBm from its harmonisable representation. Recently, Boufoussi, Dozzi and Guerbaz [13, 14] studied the existence, joint continuity and Hölder regularity of the local time of mfBm and established Chung's law of the iterated logarithm for mfBm. The main tool they applied to derive their results is the property of *one-sided* local non-determinism of mfBm.

There are multiparameter extensions of fBm, among which, two typical examples are multiparameter Lévy fBm $\xi^\alpha = \{\xi^\alpha(t), t \in \mathbb{R}^N\}$ (whose covariance function is given by (1.1) with $|\cdot|$ being the Euclidean norm in $\mathbb{R}^N$) and fractional Brownian sheets, where the former is isotropic, while the latter are anisotropic in general. Since they were introduced by Kamont [24] (see also Ayache, Léger and Pontier [4]), fractional Brownian sheets (fBs) have been studied extensively as a representative of anisotropic Gaussian random fields in recent years. See, for example, Dunker [18], Mason and Shi [32], Øksendal and Zhang [33], Xiao and Zhang [45], Ayache and Xiao [6], Ayache, Wu and Xiao [5], Wu and Xiao [39] and the references therein for further information. Still, the regularity of fBs does not evolve in the $N$-dimensional 'time' parameter $t \in \mathbb{R}_+^N$.

To model anisotropic Gaussian random fields whose regularity evolves in time, such as images, Ayache and Léger [3] and Herbin [22] introduced so-called multifractional Brownian sheets (mfBs) in terms of their moving average representations and harmonisable representations, where the constant Hurst vector of fBs is substituted by Hurst functionals. Furthermore, they showed that mfBs have a continuous modification and determined the pointwise and local Hölder exponent of mfBs. They also proved that mfBs are locally self-similar. We refer to Ayache and Léger [3] and Herbin [22] for the definitions of the corresponding concepts and results.

In studying anisotropic random fields, Xiao [44] suggested that it is more convenient to use the following metric $\rho_K$ on $\mathbb{R}^N$:

$$\rho_K(s,t) = \sum_{\ell=1}^N |s_\ell - t_\ell|^{K_\ell} \qquad \forall s,t \in \mathbb{R}^N, \tag{1.2}$$

where $K = (K_1, \ldots, K_N) \in (0,1)^N$ is a fixed vector. Denote by $H(t) = (H_1(t), \ldots, H_N(t))$ a function in $t \in \mathbb{R}_+^N$ with values in $(0,1)^N$. We say that $H(t)$ satisfies Condition A if there exist a positive number $\alpha \in (0,1)$ and a vector $(K_1, \ldots, K_N) \in (0,1)^N$ such that

A.1. for every $\ell \in \{1, \ldots, N\}$, $\alpha \leq H_\ell(t) \leq K_\ell$ for all $t \in \mathbb{R}_+^N$;
A.2. $H_\ell(t)$ ($\ell = 1, \ldots, N$) satisfies a $\rho_K$-Lipschitz condition on every compact set, that is, for every compact subset $I \subseteq \mathbb{R}_+^N$, there exist positive constants $c_\ell = c_\ell(I)$ and $\delta$ such that

$$|H_\ell(t) - H_\ell(s)| \leq c_\ell \rho_K(s,t) \qquad \forall s,t \in I \text{ with } |s-t| < \delta.$$

Now, we are ready to define multifractional Brownian sheets via their moving average representations.

***Definition 1.1.*** *Let $H(t) = (H_1(t), \ldots, H_N(t))$ be a function in $t \in \mathbb{R}_+^N$ with values in $(0,1)^N$ satisfying Condition* A. *A real-valued multifractional Brownian sheet $\{B_0^{H(t)}(t)\} = \{B_0^{H(t)}(t), t \in \mathbb{R}_+^N\}$ with functional Hurst index $H(t)$ is defined as the following moving average Wiener integral:*

$$B_0^{H(t)}(t) = \int_{\mathbb{R}^N} \prod_{\ell=1}^N [(t_\ell - u_\ell)_+^{H_\ell(t)-1/2} - (-u_\ell)_+^{H_\ell(t)-1/2}] W(\mathrm{d}u) \qquad \forall t \in \mathbb{R}_+^N, \qquad (1.3)$$

*where $s_+ = \max\{s, 0\}$, and where $W = \{W(s), s \in \mathbb{R}^N\}$ is a standard real-valued Brownian sheet.*

***Remark 1.2.*** Our Definition 1.1 generalizes the definition in Ayache and Léger [3], where they define $H_\ell(\cdot)$ as a function of $t_\ell \in \mathbb{R}_+$ ($\ell = 1, \ldots, N$). Herbin [22] defines multifractional Brownian sheets by using the following moving average representation (cf. [22] Definition 2, page 1261):

$$B_0'^{H(t)}(t) = \int_{\mathbb{R}^N} \prod_{\ell=1}^N [|t_\ell - u_\ell|^{H_\ell(t)-1/2} - |u_\ell|^{H_\ell(t)-1/2}] W(\mathrm{d}u) \qquad \forall t \in \mathbb{R}_+^N. \qquad (1.4)$$

Based on Dobrić and Ojeda [16] (see also Stoev and Taqqu [35]), we know that the multifractional Brownian sheet defined by (1.3) and the multifractional Brownian sheet defined by (1.4) have different correlation structures in general, even if $N = 1$. Our definition is more convenient to use than that of Herbin's when we derive the one-sided sectorial local non-determinism for multifractional Brownian sheets in Section 2. The form (1.3) is preferred in some applications in the one-dimensional case because it is easier to separate the future from the past.

Define the Gaussian random field $\{B^{H(t)}(t)\} = \{B^{H(t)}(t) : t \in \mathbb{R}_+^N\}$ with values in $\mathbb{R}^d$ by

$$B^{H(t)}(t) = (B_1^{H(t)}(t), \ldots, B_d^{H(t)}(t)) \qquad \forall t \in \mathbb{R}_+^N, \qquad (1.5)$$

where $\{B_1^{H(t)}(t)\}, \ldots, \{B_d^{H(t)}(t)\}$ are $d$ independent copies of $\{B_0^{H(t)}(t)\}$. Then $\{B^{H(t)}(t), t \in \mathbb{R}_+^N\}$ is called an $(N, d)$-multifractional Brownian sheet with functional Hurst index $H(t)$.

Note that if $N = 1$, then $\{B^{H(t)}(t)\}$ is a multifractional Brownian motion in $\mathbb{R}^d$ with Hurst index $H_1(t) \in (0, 1)$; if $N > 1$ and $H_1(t) \equiv H_1, \ldots, H_N(t) \equiv H_N$, then $\{B^{H(t)}(t)\}$ is an $(N, d)$-fractional Brownian sheet with Hurst index $H = (H_1, \ldots, H_N)$, which will be denoted $\{B^H(t), t \in \mathbb{R}_+^N\}$.

In this paper, we will study the existence and regularity of the local times of multifractional Brownian sheets. Our main technical tool is the property of *one-sided sectorial local* *non-determinism*. This property is more general than the sectorial local non-determinism which was first introduced by Khoshnevisan and Xiao [27] for the Brownian sheet and then extended by Wu and Xiao [39] to fractional Brownian sheets. See Section 2 for more information.

Our results show that multifractional Brownian sheets are similar to fractional Brownian sheets in many ways. They admit jointly continuous local times when the indices $H_i(t)$ stay in the range for which the fractional Brownian sheet has jointly continuous local times. We also establish a Hausdorff dimension result and, essentially, the dimension of the level set is the same as for the constant parameter case, except that we take the supremum of the constant parameter formula. We also show that the supremum can be taken locally, to establish that the fractal dimension of the random field varies in space. Hence, multifractional Brownian sheets are useful in applications such as composite materials, or porous media flow, when the material properties vary in space. They may also find useful applications in image processing.

The rest of this paper is organized as follows. In Section 2, we prove some basic results on mfBs that will be useful to our arguments. In Section 3, we provide a sufficient condition for the existence of $L^2$-local times of the $(N,d)$ mfBs and prove that the condition also implies the joint continuity of the local times. We prove the Hölder regularity of the local times in Section 4. Finally, we derive the local Hausdorff dimensions of the level sets of $\{B^{H(t)}(t)\}$ in Section 5. Our results extend the results of Ayache and Xiao [6] and Ayache, Wu and Xiao [5] for fractional Brownian sheets and Boufoussi, Dozzi and Guerbaz [13, 14] for multifractional Brownian motion to multifractional Brownian sheets.

We end the Introduction with some notation. Throughout this paper, the underlying parameter space is $\mathbb{R}^N$ or $\mathbb{R}^N_+ = [0, \infty)^N$. A parameter $t \in \mathbb{R}^N$ is written as $t = (t_1, \ldots, t_N)$, or as $\langle c \rangle$, if $t_1 = \cdots = t_N = c$. For any $s, t \in \mathbb{R}^N$ such that $s_j < t_j$ $(j = 1, \ldots, N)$ (we denote this by $s \prec t$), we define the closed interval (or rectangle) $[s, t] = \prod_{j=1}^N [s_j, t_j]$. We use $\mathcal{A}$ to represent the collection of closed intervals $[s, t]$ with $s, t \in [\varepsilon, \Lambda]^N$ for some fixed positive numbers $\varepsilon$ and $\Lambda$. For any integer $m \geq 1$, we always write $\lambda_m$ for the Lebesgue measure on $\mathbb{R}^m$ and use $\langle \cdot, \cdot \rangle$ and $|\cdot|$ to denote the ordinary scalar product and the Euclidean norm in $\mathbb{R}^m$, respectively.

Throughout this paper, an unspecified positive and finite constant will be denoted by $c$, which may not be the same in each occurrence. More specific constants in Section $i$ are numbered as $c_{i,1}, c_{i,2}, \ldots$.

## 2. Preliminaries

In this section, we will provide some lemmas that are useful for proving our main results. Lemma 2.1 is an extension of Boufoussi *et al.* [13], Lemma 3.1 from fractional Brownian motion to fractional Brownian sheets.

**Lemma 2.1.** *Let $0 < \varepsilon < \Lambda$ and $0 < \alpha < \gamma < 1$ be fixed constants. Let $\{Z_0^\kappa(t), (t, \kappa) \in \mathbb{R}_+^N \times [\alpha, \gamma]^N\}$ be a real-valued Gaussian random field defined by equation (1.3) with* $H(t) \equiv \kappa$. There then exists a constant $c_{2,1} = c(\alpha, \gamma, \varepsilon, \Lambda, N) > 0$ such that

$$\mathbb{E}[Z_0^\beta(t) - Z_0^{\beta'}(t)]^2 \leq c_{2,1} |\beta - \beta'|^2 \tag{2.1}$$

for all $t \in [\varepsilon, \Lambda]^N$ and all $\beta, \beta' \in [\alpha, \gamma]^N$.

**Proof.** For any $\beta, \beta' \in [\alpha, \gamma]^N$, we define $\kappa^0 = \beta$ and $\kappa^j = (\beta'_1, \ldots, \beta'_j, \beta_{j+1}, \ldots, \beta_N)$ for $j = 1, \ldots, N$. Clearly, $\kappa^N = \beta'$. Since

$$\mathbb{E}[Z_0^\beta(t) - Z_0^{\beta'}(t)]^2 \leq N \sum_{j=1}^{N} \mathbb{E}[Z_0^{\kappa^{j-1}}(t) - Z_0^{\kappa^j}(t)]^2, \tag{2.2}$$

it suffices for us to prove that for $j \in \{1, \ldots, N\}$ fixed,

$$\mathbb{E}[Z_0^{\kappa^{j-1}}(t) - Z_0^{\kappa^j}(t)]^2 \leq c(\alpha, \gamma, \varepsilon, \Lambda, N) |\beta_j - \beta'_j|^2. \tag{2.3}$$

By using the moving average representation of fBs (equation (1.3)), we have that

$$\mathbb{E}[Z_0^{\kappa^{j-1}}(t) - Z_0^{\kappa^j}(t)]^2$$

$$= \int_{\mathbb{R}^{N-1}} \prod_{\ell=1}^{j-1} [(t_\ell - u_\ell)_+^{\beta'_\ell - 1/2} - (-u_\ell)_+^{\beta'_\ell - 1/2}]^2$$

$$\times \prod_{\ell=j+1}^{N} [(t_\ell - u_\ell)_+^{\beta_\ell - 1/2} - (-u_\ell)_+^{\beta_\ell - 1/2}]^2 \, d\breve{u}_j \tag{2.4}$$

$$\times \int_{\mathbb{R}} [(t_j - u_j)_+^{\beta_j - 1/2} - (-u_j)_+^{\beta_j - 1/2} - ((t_j - u_j)_+^{\beta'_j - 1/2} - (-u_j)_+^{\beta'_j - 1/2})]^2 \, du_j$$

$$:= I_1 \times I_2,$$

where $\breve{u}_j = (u_1, \ldots, u_{j-1}, u_{j+1}, \ldots, u_N)$.

A change of variables shows that

$$I_1 = c \prod_{\ell=1}^{j-1} t_\ell^{2\beta_\ell} \cdot \prod_{\ell=j+1}^{N} t_\ell^{2\beta'_\ell}, \tag{2.5}$$

which is bounded for all $t \in [\varepsilon, \Lambda]^N$ and $\beta, \beta' \in [\alpha, \gamma]^N$.

Next, we estimate $I_2$. Without loss of generality, we may assume $\beta_j > \beta'_j$. We rewrite $I_2$ as the following summation:

$$I_2 = \int_0^{t_j} [(t_j - u_j)^{\beta_j - 1/2} - (t_j - u_j)^{\beta'_j - 1/2}]^2 \, du_j$$

$$+ \int_{-\infty}^{0} [((t_j - u_j)^{\beta_j - 1/2} - (-u_j)^{\beta_j - 1/2})$$

$$- ((t_j - u_j)^{\beta'_j - 1/2} - (-u_j)^{\beta'_j - 1/2})]^2 \, du_j \tag{2.6}$$

$$:= \mathrm{II}_1 + \mathrm{II}_2.$$

By the mean value theorem, we have that for some $\beta'_j \leq \eta_{j1} \leq \beta_j$ ($\eta_{j1}$ may depend on $t_j$),

$$\mathrm{II}_1 = \int_0^{t_j} [(t_j - u_j)^{\eta_{j1} - 1/2} \ln(t_j - u_j)]^2 |\beta_j - \beta'_j|^2 \, du_j$$
$$\leq c_{2,2} |\beta_j - \beta'_j|^2 \tag{2.7}$$

for all $t \in [\varepsilon, \Lambda]^N$ and $\beta, \beta' \in [\alpha, \gamma]^N$. In the above, the last inequality follows from a change of variables.

Similarly, we have that for some $\beta'_j \leq \eta_{j2} \leq \beta_j$,

$$\mathrm{II}_2 = \int_{-\infty}^0 [(t_j - u_j)^{\eta_{j2} - 1/2} \ln(t_j - u_j) - (-u_j)^{\eta_{j2} - 1/2} \ln(-u_j)]^2 |\beta_j - \beta'_j|^2 \, du_j$$
$$\leq c_{2,3} |\beta_j - \beta'_j|^2 \tag{2.8}$$

for all $t \in [\varepsilon, \Lambda]^N$ and $\beta, \beta' \in [\alpha, \gamma]^N$. Equation (2.3) is proved by combining (2.4)–(2.8). This proves Lemma 2.1. $\square$

Combining Lemma 2.1 and Lemma 8 of Ayache and Xiao [6], we have the following result.

**Lemma 2.2.** *Let $\{B_0^{H(t)}(t)\}$ be a multifractional Brownian sheet in $\mathbb{R}$. There exist positive constants $\delta > 0$, $c_{2,4}$ and $c_{2,5}$ such that for all $s, t \in [\varepsilon, \Lambda]^N$ with $|s - t| < \delta$, for any $u \in \prod_{\ell=1}^N [s_\ell \wedge t_\ell, s_\ell \vee t_\ell]$, we have*

$$c_{2,4} \sum_{\ell=1}^N |t_\ell - s_\ell|^{2H_\ell(u)} \leq \mathbb{E}[B_0^{H(t)}(t) - B_0^{H(s)}(s)]^2 \leq c_{2,5} \sum_{\ell=1}^N |t_\ell - s_\ell|^{2H_\ell(u)}. \tag{2.9}$$

**Proof.** By the elementary inequalities

$$3(a^2 + b^2 + c^2) \geq (a + b + c)^2 \geq \tfrac{1}{2} a^2 - 4b^2 - 4c^2,$$

we have that

$$3(\mathbb{E}[B_0^{H(u)}(t) - B_0^{H(u)}(s)]^2 + \mathbb{E}[B_0^{H(t)}(t) - B_0^{H(u)}(t)]^2 + \mathbb{E}[B_0^{H(s)}(s) - B_0^{H(u)}(s)]^2)$$
$$\geq \mathbb{E}[B_0^{H(t)}(t) - B_0^{H(s)}(s)]^2$$
$$\geq \tfrac{1}{2} \mathbb{E}[B_0^{H(u)}(t) - B_0^{H(u)}(s)]^2 - 4\mathbb{E}[B_0^{H(t)}(t) - B_0^{H(u)}(t)]^2$$
$$- 4\mathbb{E}[B_0^{H(u)}(s) - B_0^{H(s)}(s)]^2. \tag{2.10}$$

The first term on the right-hand side of (2.10) is the variance of the increment of a fractional Brownian sheet with Hurst index $H(u)$. By Lemma 8 in Ayache and Xiao [6], there exist positive constants $c_{2,6}$ and $c_{2,7}$, depending only on $\alpha, \gamma = \max\{K_1,\ldots,K_N\}, \varepsilon$ and $N$, such that

$$c_{2,6} \sum_{\ell=1}^{N} |t_\ell - s_\ell|^{2H_\ell(u)} \leq \mathbb{E}[B_0^{H(u)}(t) - B_0^{H(u)}(s)]^2 \leq c_{2,7} \sum_{\ell=1}^{N} |t_\ell - s_\ell|^{2H_\ell(u)}. \quad (2.11)$$

Meanwhile, by Condition A and Lemma 2.1, there exists $\delta > 0$ small such that for all $s, t \in [\varepsilon, \Lambda]^N$ with $|s - t| < \delta$, we have $|u - s| < \delta$, which implies

$$\begin{aligned}
\mathbb{E}[B_0^{H(u)}(s) - B_0^{H(s)}(s)]^2 &\leq c_{2,1} |H(u) - H(s)|^2 \\
&= c_{2,1} \sum_{\ell=1}^{N} |H_\ell(u) - H_\ell(s)|^2 \leq c_{2,8} \sum_{\ell=1}^{N} |t_\ell - s_\ell|^{2K_\ell}
\end{aligned} \quad (2.12)$$

and, similarly,

$$\mathbb{E}[B_0^{H(u)}(t) - B_0^{H(t)}(t)]^2 \leq c_{2,9} \sum_{\ell=1}^{N} |t_\ell - s_\ell|^{2K_\ell}. \quad (2.13)$$

Combining (2.10) with (2.11), (2.12) and (2.13), and noting that $\sup_{t_\ell} H_\ell(t_\ell) \leq K_\ell$ (still by Definition 1.1), we can see that there exists $\delta > 0$, which depends only on $\alpha, \gamma, \varepsilon, \Lambda$ and $N$, such that for $|s - t| < \delta$, (2.9) holds. This completes the proof of Lemma 2.2. □

**Remark 2.3.** For any fixed $t \in [\varepsilon, \Lambda]^N$ and any fixed unit vector $\eta \in \mathbb{R}^N$, let us consider the increment

$$Y_\eta(t) = B_0^{H(t+\eta)}(t + \eta) - B_0^{H(t)}(t).$$

In light of Lemma 2.2, it would be interesting to compute the bounds on the correlation function

$$r_{t,\eta}(s) = \frac{\mathbb{E}[Y_\eta(t+s)Y_\eta(t)]}{\sqrt{\mathbb{E}[Y_\eta(t+s)^2]\mathbb{E}[Y_\eta(t)^2]}} \quad (2.14)$$

and its time- and direction-varying spectral density (or Wigner–Ville distribution)

$$f_{t,\eta}(\xi) = \int_{\mathbb{R}^N} e^{-i\langle s,\xi\rangle} r_{t,\eta}(s) \, ds \quad (2.15)$$

since those are the properties of the random field that are usually used to model natural phenomena.

As we mentioned in the Introduction, the main technical tool of this paper is the property of sectorial local non-determinism. Hence, it will be helpful to say a few words about the history of various forms of local non-determinism and briefly recall some of the definitions. The concept of *local non-determinism*, which provides a powerful tool for dealing with the complex dependence structure of Gaussian processes, was first introduced by Berman [11] to unify and extend his methods for studying local times of real-valued Gaussian processes, and then extended by Pitt [34] to Gaussian random fields. The notion of *strong local non-determinism* was later introduced by Cuzick and DuPreez [15] for Gaussian processes (i.e., $N = 1$) and can be extended in various ways to Gaussian random fields (see Xiao [43]). The simplest is the following. Let $X = \{X(t), t \in \mathbb{R}^N\}$ be a real-valued Gaussian random field with $0 < \mathbb{E}[X(t)^2] < \infty$ for all $t \in I$, where $I \in \mathcal{A}$ is an interval. Let $\phi$ be a given function such that $\phi(0) = 0$ and $\phi(r) > 0$ for $r > 0$. Then $X$ is said to be *strongly locally $\phi$-non-deterministic* (SL$\phi$ND) on $I$ if there exist positive constants $c_{2,10}$ and $r_0$ such that for all $t \in I$ and all $0 < r \leq \min\{|t|, r_0\}$,

$$\mathrm{Var}(X(t)|X(s) : s \in I, r \leq |s - t| \leq r_0) \geq c_{2,10} \phi(r). \tag{2.16}$$

Here, $\mathrm{Var}(X(t)|X(s), s \in S)$ denotes the conditional variance of $X(t)$, given $X(s), s \in S$. Pitt [34] proved that fractional Brownian motion $\xi^\alpha$ is SL$\phi$ND with $\phi(r) = r^{2\alpha}$. When $N = 1$ and $\mathrm{Var}(X(t)|X(s) : s \in I, r \leq t - s \leq r_0) \geq c_{2,10} \phi(r)$, $X$ is said to satisfy *one-sided strong local $\phi$-non-determinism*. The properties of strong local non-determinism have played important roles in studying the regularity of local times, small ball probabilities, uniform Hausdorff dimension results and other sample path properties of Gaussian processes and Gaussian random fields. We refer to Xiao [43, 44] for more information on applications of SL$\phi$ND, as well as other forms of local non-determinism.

It is well known that the Brownian sheet does not satisfy the property of local non-determinism in the sense of Berman and Pitt (hence, it does not satisfy (2.16) either). Recently, Khoshnevisan and Xiao [27] introduced the concept of *sectorial local non-determinism* and proved that the Brownian sheet has the property of sectorial local non-determinism. Wu and Xiao [39] extended the result of Khoshnevisan and Xiao [27] and proved that the fractional Brownian sheet $\{B_0^H(t), t \in \mathbb{R}_+^N\}$, where $H = (H_1, \ldots, H_N) \in (0, 1)$, also has the property of sectorial local non-determinism. Namely, for any fixed positive number $\varepsilon \in (0, 1)$, there exists a positive constant $c_{2,11}$, depending only on $\varepsilon, H$ and $N$, such that for all positive integers $n \geq 2$, and all $t^1, \ldots, t^n \in [\varepsilon, \infty)^N$, we have

$$\mathrm{Var}(B_0^H(t^n)|B_0^H(t^1), \ldots, B_0^H(t^{n-1})) \geq c_{2,11} \sum_{j=1}^N \min_{0 \leq k \leq n-1} |t_j^n - t_j^k|^{2H_j}, \tag{2.17}$$

where $t^0 = 0$.

The concept of sectorial local non-determinism opens the door to unifying the previously different treatments for the Brownian sheet and fractional Brownian motion, and can also be used to prove new results about fractional Brownian motion. See, for example, Khoshnevisan, Wu and Xiao [26], Wu and Xiao [39] and Ayache, Wu and Xiao [5] for studies on sample path properties of the Brownian sheet and fractional Brownian sheets by using the sectorial local non-determinism property.

It is not known whether the multifractional Brownian sheet $\{B_0^{H(t)}(t)\}$ is sectorially locally non-deterministic. To study the local times of $(N,d)$-mfBs $\{B^{H(t)}(t)\}$, we will work with multifractional Liouville sheets (mfLs) at first and prove that they can be well approximated by Gaussian processes which have the property of one-sided strong local non-determinism (see (2.22) and (2.23) below). For a function $H(t) = (H_1(t), \ldots, H_N(t))$ satisfying Condition A, the real-valued, centered Gaussian random field $\{X_0^{H(t)}(t)\} = \{X_0^{H(t)}(t), t \in \mathbb{R}_+^N\}$ defined by

$$X_0^{H(t)}(t) = \int_{[0,t]} \prod_{\ell=1}^N (t_\ell - s_\ell)^{H_\ell(t) - 1/2} W(\mathrm{d}s), \qquad t \in \mathbb{R}_+^N, \tag{2.18}$$

is called a multifractional Liouville sheet with functional Hurst index $H(t)$. One parameter mfLs were first introduced by Lim and Muniandy [30] as an extension of fBm; see Lim [29] for more properties on one parameter mfLs.

It follows from (1.3) that for every $t \in \mathbb{R}_+^N$,

$$B_0^{H(t)}(t) = X_0^{H(t)}(t) + \int_{(-\infty,t]\setminus[0,t]} \prod_{\ell=1}^N g_\ell(t_\ell, s_\ell) W(\mathrm{d}s), \tag{2.19}$$

where

$$g_\ell(t_\ell, s_\ell) = ((t_\ell - s_\ell)_+)^{H_\ell(t) - 1/2} - ((-s_\ell)_+)^{H_\ell(t) - 1/2}$$

and the two fields on the right-hand side of (2.19) are independent. We will show that in studying the regularity properties of the local times of $\{B_0^{H(t)}(t)\}$, the Liouville sheet $\{X_0^{H(t)}(t)\}$ plays a crucial role and the second field in (2.19) can be neglected. More precisely, we will make use of the following property. For all integers $n \geq 2$, $t^1, \ldots, t^n \in \mathbb{R}_+^N$ and $u_1, \ldots, u_n \in \mathbb{R}$, we have

$$\mathrm{Var}\left(\sum_{j=1}^n u_j B_0^{H(t^j)}(t^j)\right) \geq \mathrm{Var}\left(\sum_{j=1}^n u_j X_0^{H(t^j)}(t^j)\right). \tag{2.20}$$

Next, we use an argument in Ayache and Xiao [6] to provide a useful decomposition for $\{X_0^{H(t)}(t)\}$. For every $t \in [\varepsilon, \Lambda]^N$, we decompose the rectangle $[0,t]$ into the following disjoint union of subrectangles:

$$[0,t] = [0,\varepsilon]^N \cup \bigcup_{\ell=1}^N R_\ell(t) \cup \Delta(\varepsilon, t), \tag{2.21}$$

where $R_\ell(t) = \{r \in [0,\Lambda]^N : 0 \leq r_i \leq \varepsilon \text{ if } i \neq \ell, \varepsilon < r_\ell \leq t_\ell\}$ and $\Delta(\varepsilon, t)$ can be written as a union of $2^N - N - 1$ subrectangles of $[0,t]$. Denote the integrand in (2.18) by $g(t,r)$. It follows from (2.21) that for every $t \in [\varepsilon, \infty)^N$,

$$X_0^{H(t)}(t) = \int_{[0,\varepsilon]^N} g(t,r) W(\mathrm{d}r) + \sum_{\ell=1}^N \int_{R_\ell(t)} g(t,r) W(\mathrm{d}r) + \int_{\Delta(\varepsilon,t)} g(t,r) W(\mathrm{d}r)$$

$$:= X(\varepsilon,t) + \sum_{\ell=1}^N Y_\ell(t) + Z(\varepsilon,t). \tag{2.22}$$

Since the processes $X(\varepsilon,t)$, $Y_\ell(t)$ $(1 \leq \ell \leq N)$ and $Z(\varepsilon,t)$ are defined by the stochastic integrals with respect to $W$ over disjoint sets, they are independent Gaussian random fields.

The following lemma shows that every process $Y_\ell(t)$ has the property of strong local non-determinism along the $\ell$th direction. It will be essential to our proofs.

**Lemma 2.4.** *Let $I \in \mathcal{A}$ and let $\ell \in \{1, 2, \ldots, N\}$ be fixed. For all integers $n \geq 2$ and $t^1, \ldots, t^n \in I$ such that*

$$t_\ell^1 \leq t_\ell^2 \leq \cdots \leq t_\ell^n,$$

*we have*

$$\mathrm{Var}(Y_\ell(t^n) | Y_\ell(t^j) : 0 \leq j \leq n-1) \geq c_{2,12} |t_\ell^n - t_\ell^{n-1}|^{2H_\ell(t^n)}, \tag{2.23}$$

*where $t_\ell^0 = 0$ and $c_{2,12} > 0$ is a constant depending only on $\varepsilon$ and $I$.*

**Proof.** The proof is in the same spirit as the proof of Lemma 2.1 of Ayache, Wu and Xiao [5]. Working in the Hilbert space setting, the conditional variance in (2.23) is the square of the $L^2(\mathbb{P})$-distance of $Y_\ell(t^n)$ from the subspace generated by $Y_\ell(t^j)(0 \leq j \leq n-1)$. Hence, it is sufficient to show that there exists a constant $c_{2,12}$ such that

$$\mathbb{E}\left( Y_\ell(t^n) - \sum_{j=1}^{n-1} a_j Y_\ell(t^j) \right)^2 \geq c_{2,12} |t_\ell^n - t_\ell^{n-1}|^{2H_\ell(t^n)} \tag{2.24}$$

for all $a_j \in \mathbb{R}$ $(j = 1, \ldots, n-1)$. However, by splitting $R_\ell(t^n)$ into two disjoint parts and using the independence, we derive that

$$\mathbb{E}\left( Y_\ell(t^n) - \sum_{j=0}^{n-1} a_j Y_\ell(t^j) \right)^2 \geq \mathbb{E}\left( \int_{R_\ell(t^n) \setminus R_\ell(t^{n-1})} g(t^n,r) W(\mathrm{d}r) \right)^2$$

$$\geq \int_0^\varepsilon \cdots \int_{t_\ell^{n-1}}^{t_\ell^n} \cdots \int_0^\varepsilon \prod_{k=1}^N (t_k^n - r_k)^{2H_k(t^n)-1} \, \mathrm{d}r \tag{2.25}$$

$$\geq c_{2,12} |t_\ell^n - t_\ell^{n-1}|^{2H_\ell(t^n)}.$$

This proves (2.24) and hence Lemma 2.4. □

Combining Lemma 2.2 and Lemma 2.4 with the proofs of Lemma 2.1 and Lemma 8.1 in Berman [11] (see also Boufoussi, Dozzi and Guerbaz [13], Theorem 3.3), we have the following.

**Proposition 2.5.** *For every integer $n \geq 2$, there exist positive constants $C_n$ and $\delta$ (both of these may depend on $n$) such that for every $\ell = 1, \ldots, N$,*

$$\mathrm{Var}\left(\sum_{j=1}^{n} u_j [Y_\ell(t^j) - Y_\ell(t^{j-1})]\right) \geq C_n \sum_{j=1}^{n} u_j^2 \mathrm{Var}[Y_\ell(t^j) - Y_\ell(t^{j-1})] \qquad (2.26)$$

*for all $(u_1, \ldots, u_n) \in \mathbb{R}^n$ and all points $t^1, \ldots, t^n \in I$ satisfying $t_\ell^1 < \cdots < t_\ell^n$ with $t_\ell^n - t_\ell^1 < \delta$.*

The following lemma relates the multifractional Brownian sheet $\{B_0^{H(t)}(t)\}$ to the independent Gaussian random fields $Y_\ell$ ($\ell = 1, \ldots, N$), which is a direct extension of Ayache, Wu and Xiao [5], Lemma 2.2.

**Lemma 2.6.** *Let $I \in \mathcal{A}$. For all integers $n \geq 2$, $t^1, \ldots, t^n \in I$ and $u_1, \ldots, u_n \in \mathbb{R}$, we have*

$$\mathrm{Var}\left(\sum_{j=1}^{n} u_j B_0^{H(t^j)}(t^j)\right) \geq \sum_{\ell=1}^{N} \mathrm{Var}\left(\sum_{j=1}^{n} u_j Y_\ell(t^j)\right). \qquad (2.27)$$

*Consequently, for all $k \in \{1, \ldots, N\}$ and positive numbers $p_1, \ldots, p_k \geq 1$ satisfying $\sum_{\ell=1}^{k} p_\ell^{-1} = 1$, we have*

$$\frac{1}{[\det \mathrm{Cov}(B_0^{H(t^1)}(t^1), \ldots, B_0^{H(t^n)}(t^n))]^{1/2}} \leq \prod_{\ell=1}^{k} \frac{c_{2,13}^n}{[\det \mathrm{Cov}(Y_\ell(t^1), \ldots, Y_\ell(t^n))]^{1/(2p_\ell)}}, \qquad (2.28)$$

*where $\det \mathrm{Cov}(Z_1, \ldots, Z_n)$ denotes the determinant of the covariance matrix of the Gaussian random vector $(Z_1, \ldots, Z_n)$.*

**Remark 2.7.** By (2.27) and Lemma 2.4, we derive that $\{B_0^{H(t)}(t)\}$ satisfies the following property: for all integers $n \geq 2$ and $t^1, \ldots, t^n \in I$ such that $t^j \prec t^n$ for every $j \leq n-1$,

$$\mathrm{Var}(B_0^H(t^n) | B_0^H(t^1), \ldots, B_0^H(t^{n-1})) \geq c_{2,12} \sum_{j=1}^{N} \min_{0 \leq k \leq n-1}(t_j^n - t_j^k)^{2H_j(t^n)}. \qquad (2.29)$$

This property is weaker than (2.17) and will be referred to as *one-sided sectorial local non-determinism*.

We will also make use of the following technical lemmas, among which Lemma 2.8 is from Xiao and Zhang [45], Lemma 2.9 is proved in Ayache and Xiao [6] and Lemma 2.10 and Lemma 2.11 are from Ayache, Wu and Xiao [5].

**Lemma 2.8.** *Let $0 < h < 1$ be a constant. Then for any $\delta > 2h$, $M > 0$ and $\beta > 0$, there exists a positive and finite constant $c_{2,14}$, depending only on $\delta, \varepsilon, \beta$ and $M$, such that for all $0 < a \leq M$,*

$$\int_{\varepsilon}^{1} dr \int_{\varepsilon}^{1} [a + |s - r|^{2h}]^{-\beta} \, ds \leq c_{2,14}(a^{-(\beta - 1/\delta)} + 1). \tag{2.30}$$

**Lemma 2.9.** *Let $\alpha, \beta$ and $\eta$ be positive constants. For $A > 0$ and $B > 0$, let*

$$J := J(A, B) = \int_{0}^{1} \frac{dt}{(A + t^{\alpha})^{\beta}(B + t)^{\eta}}. \tag{2.31}$$

*There then exist finite constants $c_{2,15}$ and $c_{2,16}$, depending only on $\alpha, \beta$ and $\eta$, such that the following hold for all reals $A, B > 0$ satisfying $A^{1/\alpha} \leq c_{2,15} B$:*

(i) *if $\alpha\beta > 1$, then*

$$J \leq c_{2,16} \frac{1}{A^{\beta - \alpha^{-1}} B^{\eta}}; \tag{2.32}$$

(ii) *if $\alpha\beta = 1$, then*

$$J \leq c_{2,16} \frac{1}{B^{\eta}} \log(1 + BA^{-1/\alpha}); \tag{2.33}$$

(iii) *if $0 < \alpha\beta < 1$ and $\alpha\beta + \eta \neq 1$, then*

$$J \leq c_{2,16} \left( \frac{1}{B^{\alpha\beta + \eta - 1}} + 1 \right). \tag{2.34}$$

**Lemma 2.10.** *Let $(\rho_1, \ldots, \rho_N) \subseteq (0, 1)^N$. For any $q \in [0, \sum_{\ell=1}^{N} \rho_\ell^{-1})$, let $\tau \in \{1, \ldots, N\}$ be the integer such that*

$$\sum_{\ell=1}^{\tau - 1} \frac{1}{\rho_\ell} \leq q < \sum_{\ell=1}^{\tau} \frac{1}{\rho_\ell}, \tag{2.35}$$

*with the convention that $\sum_{\ell=1}^{0} \frac{1}{\rho_\ell} := 0$. There then exists a positive constant $\Delta_\tau \leq 1$, depending only on $(\rho_1, \ldots, \rho_N)$, such that for every $\Delta \in (0, \Delta_\tau)$, we can find $\tau$ real numbers $p_\ell \geq 1$ $(1 \leq \ell \leq \tau)$ satisfying the following properties:*

$$\sum_{\ell=1}^{\tau} \frac{1}{p_\ell} = 1, \qquad \frac{\rho_\ell q}{p_\ell} < 1 \qquad \forall \ell = 1, \ldots, \tau \tag{2.36}$$

*and*

$$(1 - \Delta) \sum_{\ell=1}^{\tau} \frac{\rho_\ell q}{p_\ell} \leq \rho_\tau q + \tau - \sum_{\ell=1}^{\tau} \frac{\rho_\tau}{\rho_\ell}. \tag{2.37}$$

Furthermore, if we let $\alpha_\tau := \sum_{\ell=1}^{\tau} \frac{1}{\rho_\ell} - q > 0$, then for any positive number $\rho \in (0, \frac{\alpha_\tau}{2\tau})$, there exists an $\ell_0 \in \{1, \ldots, \tau\}$ such that

$$\frac{\rho_{\ell_0} q}{p_{\ell_0}} + 2\rho_{\ell_0}\rho < 1. \tag{2.38}$$

**Lemma 2.11.** *For all integers $n \geq 1$, positive numbers $a$, $r$, $0 < b_j < 1$ and an arbitrary $s_0 \in [0, a/2]$,*

$$\int_{a \leq s_1 \leq \cdots \leq s_n \leq a+r} \prod_{j=1}^{n}(s_j - s_{j-1})^{-b_j}\,\mathrm{d}s_1 \cdots \mathrm{d}s_n \leq c_{2,17}^n (n!)^{(1/n)\sum_{j=1}^{n} b_j - 1} r^{n - \sum_{j=2}^{n} b_j}, \tag{2.39}$$

*where $c_{2,17} > 0$ is a constant depending only on $a$ and the $b_j$'s. In particular, if $b_j = \alpha$ for all $j = 1, \ldots, n$, then*

$$\int_{a \leq s_1 \leq \cdots \leq s_n \leq a+r} \prod_{j=1}^{n}(s_j - s_{j-1})^{-\alpha}\,\mathrm{d}s_1 \cdots \mathrm{d}s_n \leq c_{2,17}^n (n!)^{\alpha - 1} r^{n(1-(1-1/n)\alpha)}. \tag{2.40}$$

Finally, we conclude this section by briefly recalling some aspects of the theory of local times. For excellent surveys on local times of random and deterministic vector fields, we refer to Geman and Horowitz [20] and Dozzi [17].

Let $X(t)$ be a Borel vector field on $\mathbb{R}^N$ with values in $\mathbb{R}^d$. For any Borel set $I \subseteq \mathbb{R}^N$, the *occupation measure* of $X$ on $I$ is defined as the following measure on $\mathbb{R}^d$:

$$\mu_I(\bullet) = \lambda_N\{t \in I : X(t) \in \bullet\}.$$

If $\mu_I$ is absolutely continuous with respect to $\lambda_d$, we say that $X(t)$ has *local times* on $I$ and define its local times, $L(\bullet, I)$, as the Radon–Nikodým derivative of $\mu_I$ with respect to $\lambda_d$, that is,

$$L(x, I) = \frac{\mathrm{d}\mu_I}{\mathrm{d}\lambda_d}(x) \qquad \forall x \in \mathbb{R}^d.$$

In the above, $x$ is the so-called *space variable* and $I$ is the *time variable*. Sometimes, we write $L(x, t)$ in place of $L(x, [0, t])$. Note that if $X$ has local times on $I$, then for every Borel set $J \subseteq I$, $L(x, J)$ also exists.

By standard martingale and monotone class arguments, one can deduce that the local times have a measurable modification that satisfies the following *occupation density formula* (see Geman and Horowitz [20], Theorem 6.4): for every Borel set $I \subseteq \mathbb{R}^N$ and every measurable function $f : \mathbb{R}^d \to \mathbb{R}_+$,

$$\int_I f(X(t))\,\mathrm{d}t = \int_{\mathbb{R}^d} f(x)L(x, I)\,\mathrm{d}x. \tag{2.41}$$

Suppose we fix a rectangle $I = \prod_{i=1}^{N}[a_i, a_i + h_i]$. Then, whenever we can choose a version of the local time, still denoted by $L(x, \prod_{i=1}^{N}[a_i, a_i + t_i])$, such that it is a continuous function of $(x, t_1, \ldots, t_N) \in \mathbb{R}^d \times \prod_{i=1}^N [0, h_i]$, $X$ is said to have a *jointly continuous local time* on $I$. When a local time is jointly continuous, $L(x, \bullet)$ can be extended to be a finite Borel measure supported on the level set

$$X_I^{-1}(x) = \{t \in I : X(t) = x\}; \tag{2.42}$$

see Adler [1] for details. In other words, local times often act as a natural measure on the level sets of $X$. Hence, they are useful in studying the various fractal properties of level sets and inverse images of the vector field $X$. In this regard, we refer to Berman [10], Ehm [19] and Xiao [41].

It follows from Geman and Horowitz [20], ((25.5) and (25.7)) (see also Pitt [34]) that for all $x, y \in \mathbb{R}^d$, $I \in \mathcal{A}$ and all integers $n \geq 1$,

$$\mathbb{E}[L(x, I)^n] = (2\pi)^{-nd} \int_{I^n} \int_{\mathbb{R}^{nd}} \exp\left(-\mathrm{i} \sum_{j=1}^n \langle u^j, x \rangle\right)$$
$$\times \mathbb{E} \exp\left(\mathrm{i} \sum_{j=1}^n \langle u^j, X(t^j) \rangle\right) \mathrm{d}\overline{u}\, \mathrm{d}\overline{t} \tag{2.43}$$

and for all even integers $n \geq 2$,

$$\mathbb{E}[(L(x, I) - L(y, I))^n] = (2\pi)^{-nd} \int_{I^n} \int_{\mathbb{R}^{nd}} \prod_{j=1}^n \left[\mathrm{e}^{-\mathrm{i}\langle u^j, x\rangle} - \mathrm{e}^{-\mathrm{i}\langle u^j, y\rangle}\right]$$
$$\times \mathbb{E} \exp\left(\mathrm{i} \sum_{j=1}^n \langle u^j, X(t^j) \rangle\right) \mathrm{d}\overline{u}\, \mathrm{d}\overline{t}, \tag{2.44}$$

where $\overline{u} = (u^1, \ldots, u^n), \overline{t} = (t^1, \ldots, t^n)$ and each $u^j \in \mathbb{R}^d, t^j \in I \subseteq (0, \infty)^N$. In the coordinate notation, we then write $u^j = (u_1^j, \ldots, u_d^j)$.

## 3. Local times: existence and joint continuity

In this section, we consider the existence and regularity of local times of mfBs. Theorem 3.1 and Corollary 3.2 provide sufficient conditions for the existence of local times. Then in Theorem 3.4 we prove that, under the same condition as in Theorem 3.1, the local times of mfBs have a jointly continuous version.

For $I \in \mathcal{A}$, let $\overline{H} = (\overline{H}_1, \ldots, \overline{H}_N)$ be the vector defined by

$$\overline{H}_\ell = \max_{t \in I} H_\ell(t) \qquad \text{for } \ell = 1, \ldots, N. \tag{3.1}$$

The index $\overline{H}$ depends on $I$, but, for simplicity, we have deleted $I$ from the notation.

**Theorem 3.1.** *Let $I \in \mathcal{A}$ and $\overline{H}$ be the vector defined in (3.1). If $d < \sum_{\ell=1}^{N} \frac{1}{\overline{H}_\ell}$, then mfBs $\{B^{H(t)}(t)\}$ admits an $L^2$-integrable local time $L(\cdot, I)$ almost surely.*

**Proof.** Without loss of generality, we may, and will, assume that $\delta = \text{diam}(I)$ is sufficiently small that (2.9) holds for all $s, t \in I$. In particular, we assume $\delta \leq 1$.

To prove the existence of local time on $I$, by Geman and Horowitz [20] Theorem 21.9, it suffices to prove that

$$\int_I \int_I (\mathbb{E}[B_0^{H(t)}(t) - B_0^{H(s)}(s)]^2)^{-d/2} \, ds \, dt < \infty. \tag{3.2}$$

It follows from (2.9) and a change of variables that

$$\int_I \int_I (\mathbb{E}[B_0^{H(t)}(t) - B_0^{H(s)}(s)]^2)^{-d/2} \, ds \, dt \leq c \int_I \int_I \left( \sum_{\ell=1}^{N} |t_\ell - s_\ell|^{2H_\ell(s)} \right)^{-d/2} ds \, dt$$

$$\leq c_{3,1} \int_I \int_I \left( \sum_{\ell=1}^{N} |t_\ell - s_\ell|^{2\overline{H}_\ell} \right)^{-d/2} ds \, dt. \tag{3.3}$$

Since $d < \sum_{\ell=1}^{N} \frac{1}{\overline{H}_\ell}$, we can estimate the last integral in a way similar to the proof of Xiao and Zhang [45] Theorem 3.6. Namely, by applying Lemma 2.8 repeatedly, we derive that the last integral is finite. This proves (3.2) and hence Theorem 3.1. $\square$

The following is a consequence of Theorem 3.1 which gives a more natural condition for the existence of local times of mfBs.

**Corollary 3.2.** *Let $I \in \mathcal{A}$ be fixed. If $\sum_{\ell=1}^{N} \frac{1}{H_\ell(t)} > d$ for all $t \in I$, then mfBs $\{B^{H(t)}(t)\}$ admits an $L^2$-integrable local time $L(\cdot, I)$ almost surely.*

**Proof.** Since the functions $H_1(t), \ldots, H_N(t)$ are uniformly continuous on $I$, we can divide $I$ into subintervals $\{I_p\}$ such that for each $I_p$, we have $\sum_{\ell=1}^{N} \frac{1}{\overline{H}_\ell(p)} > d$, where $\overline{H}_\ell(p) = \max_{t \in I_p} H_\ell(t)$. It follows from Theorem 3.1 that on every $I_p$, mfBs $\{B^{H(t)}(t)\}$ has an $L^2$-integrable local time $L(\cdot, I_p)$ almost surely. This implies that $\{B^{H(t)}(t)\}$ has an $L^2$-integrable local time on $I$, which concludes the proof of Corollary 3.2. $\square$

*Remark 3.3.* The condition in Corollary 3.2 is almost the best possible, in the sense that, if $d > \sum_{\ell=1}^{N} \frac{1}{H_\ell(t)}$ for some $t \in I$, then it can be proven using Geman and Horowitz [20] Theorem 21.9 that $\{B^{H(t)}(t)\}$ has no $L^2(\mathbb{R}^d \times \Omega)$-integrable local times on $I$. In the case when $d \leq \sum_{\ell=1}^{N} \frac{1}{H_\ell(t)}$ for all $t \in I$, but the equality only holds for $t$ in a set of Lebesgue measure 0, the existence of local times is rather subtle and requires the imposition of further assumptions on $(H_1(t), \ldots, H_N(t))$. Hence, it will not be discussed here.

We now consider the joint continuity of the local times of $\{B^{H(t)}(t)\}$. For convenience, we first prove that under the same condition as in Theorem 3.1, the local time of $\{B^{H(t)}(t)\}$ has a version that is jointly continuous in both space and time variables. We then apply the same argument as in the proof of Corollary 3.2 to show that the same conclusion holds, provided $\sum_{\ell=1}^{N} \frac{1}{H_\ell(t)} > d$ for all $t \in I$. Our results extend those of Ehm [19] for the Brownian sheet and of Ayache, Wu and Xiao [5] for fractional Brownian sheets.

**Theorem 3.4.** *Let $\{B^{H(t)}(t)\} = \{B^{H(t)}(t), t \in \mathbb{R}_+^N\}$ be a multifractional Brownian sheet with values in $\mathbb{R}^d$. Let $I \in \mathcal{A}$ and let $\overline{H}$ be the vector defined in (3.1). If $d < \sum_{\ell=1}^{N} \frac{1}{\overline{H}_\ell}$, then $\{B^{H(t)}(t)\}$ has a jointly continuous local time on $I$.*

The main idea for proving Theorem 3.4 is similar to those in Ehm [19], Xiao [41] and Ayache, Wu and Xiao [5]. That is, we first apply the Fourier analytic arguments to derive estimates for the moments of the local times and then apply a multiparameter version of the Kolmogorov continuity theorem (cf. Khoshnevisan [25]). As in Ayache, Wu and Xiao [5], the 'one-sided' sectorial local non-deterministic properties of multifractional Liouville sheets proved in Section 2 (see Lemma 2.4 and Proposition 2.5) will play important roles in deriving moment estimates in Lemmas 3.5 and 3.7 below. However, due to the non-stationarity and the lack of *two-sided* local non-determinism of multifractional Brownian sheets, we need to make several modifications in our proofs.

For convenience, we further assume

$$0 < \overline{H}_1 \leq \cdots \leq \overline{H}_N < 1. \tag{3.4}$$

**Lemma 3.5.** *Assume the conditions of Theorem 3.4 hold. Let $\tau$ be the unique integer in $\{1,\ldots,N\}$ satisfying*

$$\sum_{\ell=1}^{\tau-1} \frac{1}{\overline{H}_\ell} \leq d < \sum_{\ell=1}^{\tau} \frac{1}{\overline{H}_\ell}. \tag{3.5}$$

*There then exists a positive constant $c_{3,2}$, depending only on $N$, $d$, $\overline{H}$ and $I$, such that for all $x \in \mathbb{R}^d$, all subintervals $T = [a, a + \langle h \rangle] \subseteq I$ with $h > 0$ small and all integers $n \geq 1$,*

$$\mathbb{E}[L(x,T)^n] \leq c_{3,2}^n (n!)^{N-\beta_\tau} h^{n\beta_\tau}, \tag{3.6}$$

*where $\beta_\tau = N - \tau - \overline{H}_\tau d + \sum_{\ell=1}^{\tau} \overline{H}_\tau / \overline{H}_\ell$.*

**Proof.** For later use, we will start with an arbitrary closed interval $T = \prod_{\ell=1}^{N} [a_\ell, a_\ell + h_\ell] \subseteq I$. It follows from (2.43) and the fact that $\{B_1^{H(t)}(t)\},\ldots,\{B_d^{H(t)}(t)\}$ are independent copies of $\{B_0^{H(t)}(t)\}$ that for all integers $n \geq 1$,

$$\mathbb{E}[L(x,T)^n] \leq (2\pi)^{-nd} \int_{T^n} \prod_{k=1}^{d} \left\{ \int_{\mathbb{R}^n} \exp\left[-\tfrac{1}{2} \operatorname{Var}\left(\sum_{j=1}^{n} u_k^j B_0^{H(t^j)}(t^j)\right)\right] \mathrm{d}U_k \right\} \mathrm{d}\bar{t}, \tag{3.7}$$

where $U_k = (u_k^1, \ldots, u_k^n) \in \mathbb{R}^n$. Fix $k = 1, \ldots, d$ and denote the inner integral in (3.7) by $\mathcal{N}_k$. Then by Lemma 2.6, we have

$$\begin{aligned}\mathcal{N}_k &\leq \int_{\mathbb{R}^n} \exp\left[-\tfrac{1}{2}\sum_{\ell=1}^{N} \mathrm{Var}\left(\sum_{j=1}^{n} u_k^j Y_\ell(t^j)\right)\right] \mathrm{d}U_k \\ &\leq \int_{\mathbb{R}^n} \exp\left[-\tfrac{1}{2}\sum_{\ell=1}^{\tau} \mathrm{Var}\left(\sum_{j=1}^{n} u_k^j Y_\ell(t^j)\right)\right] \mathrm{d}U_k.\end{aligned} \quad (3.8)$$

Since (3.5) holds, we apply Lemma 2.10 with $\Delta = n^{-1}$ and $q = d$ to obtain $\tau$ positive numbers $p_1, \ldots, p_\tau \geq 1$ satisfying (2.36) and (2.37).

Applying the generalized Hölder inequality (Hardy [21] page 140), we derive that

$$\begin{aligned}\mathcal{N}_k &\leq \prod_{\ell=1}^{\tau}\left\{\int_{\mathbb{R}^n} \exp\left[-\frac{p_\ell}{2} \mathrm{Var}\left(\sum_{j=1}^{n} u_k^j Y_\ell(t^j)\right)\right] \mathrm{d}U_k\right\}^{1/p_\ell} \\ &= c_{3,3}^n \prod_{\ell=1}^{\tau} [\det \mathrm{Cov}(Y_\ell(t^1), \ldots, Y_\ell(t^n))]^{-1/(2p_\ell)},\end{aligned} \quad (3.9)$$

where the last equality follows from the fact that $(Y_\ell(t^1), \ldots, Y_\ell(t^n))$ is a Gaussian vector with mean 0. Hence, it follows from (3.7) and (3.9) that

$$\mathbb{E}[L(x,T)^n] \leq c_{3,3}^n \int_{T^n} \prod_{\ell=1}^{\tau} [\det \mathrm{Cov}(Y_\ell(t^1), \ldots, Y_\ell(t^n))]^{-d/(2p_\ell)} \, \mathrm{d}\bar{t}. \quad (3.10)$$

To evaluate the integral in (3.10), we will first integrate $[\mathrm{d}t_\ell^1 \cdots \mathrm{d}t_\ell^n]$ for $\ell = 1, \ldots, \tau$. To this end, we use the following fact about multivariate normal distributions: for any Gaussian random vector $(Z_1, \ldots, Z_n)$,

$$\det \mathrm{Cov}(Z_1, \ldots, Z_n) = \mathrm{Var}(Z_1) \prod_{j=2}^{n} \mathrm{Var}(Z_j | Z_1, \ldots, Z_{j-1}). \quad (3.11)$$

By the above fact and Lemma 2.4, we can derive that for every $\ell \in \{1, \ldots, \tau\}$ and for all $t^1, \ldots, t^n \in T = \prod_{\ell=1}^{N}[a_\ell, a_\ell + h_\ell]$ satisfying

$$a_\ell \leq t_\ell^{\pi_\ell(1)} \leq t_\ell^{\pi_\ell(2)} \leq \cdots \leq t_\ell^{\pi_\ell(n)} \leq a_\ell + h_\ell \quad (3.12)$$

for some permutation $\pi_\ell$ of $\{1, \ldots, N\}$, we have

$$\det \mathrm{Cov}(Y_\ell(t^1), \ldots, Y_\ell(t^n)) \geq c_{3,4}^n \prod_{j=1}^{n} (t_\ell^{\pi_\ell(j)} - t_\ell^{\pi_\ell(j-1)})^{2H_\ell(t^{\pi_\ell(j)})} \quad (3.13)$$

$$\geq c_{3,5}^n \prod_{j=1}^n (t_\ell^{\pi_\ell(j)} - t_\ell^{\pi_\ell(j-1)})^{2\overline{H}_\ell},$$

where $t_\ell^{\pi_\ell(0)} := \varepsilon$ (recall the decomposition (2.22)).

We choose $\varepsilon < \frac{1}{2}\min\{a_\ell, 1 \leq \ell \leq N\}$ so that Lemma 2.11 is applicable. It follows from (3.12) and (3.13) that

$$\int_{[a_\ell, a_\ell+h_\ell]^n} [\det \mathrm{Cov}(Y_\ell(t^1), \ldots, Y_\ell(t^n))]^{-d/(2p_\ell)} \, \mathrm{d}t_\ell^1 \cdots \mathrm{d}t_\ell^n$$
$$\leq \sum_{\pi_\ell} c^n \int_{a_\ell \leq t_\ell^{\pi_\ell(1)} \leq \cdots \leq t_\ell^{\pi_\ell(n)} \leq a_\ell+h_\ell} \prod_{j=1}^n \frac{1}{(t_\ell^{\pi_\ell(j)} - t_\ell^{\pi_\ell(j-1)})^{\overline{H}_\ell d/p_\ell}} \, \mathrm{d}t_\ell^1 \cdots \mathrm{d}t_\ell^n \quad (3.14)$$
$$\leq c_{3,6}^n (n!)^{\overline{H}_\ell d/p_\ell} h_\ell^{n(1-(1-1/n)\overline{H}_\ell d/p_\ell)}.$$

In the above, the last inequality follows from (2.40) (recall that $\overline{H}_\ell d/p_\ell < 1$).

Combining (3.10) and (3.14), we have

$$\mathbb{E}[L(x,T)^n] \leq c_{3,7}^n (n!)^{\sum_{\ell=1}^\tau \overline{H}_\ell d/p_\ell} \prod_{\ell=1}^\tau h_\ell^{n(1-(1-1/n)\overline{H}_\ell d/p_\ell)} \cdot \prod_{\ell=\tau+1}^N h_\ell^n. \quad (3.15)$$

We now consider the special case when $T = [a, a + \langle h \rangle]$, that is, $h_1 = \cdots = h_N = h$ with $h < \delta$, (3.15) and (2.37) with $\Delta = n^{-1}$ and $q = d$ together yield

$$\mathbb{E}[L(x,T)^n] \leq c_{3,8}^n (n!)^{\sum_{\ell=1}^\tau \overline{H}_\ell d/p_\ell} h^{n(N-(1-n^{-1})\sum_{\ell=1}^\tau \overline{H}_\ell d/p_\ell)}$$
$$\leq c_{3,9}^n (n!)^{N-\beta_\tau} h^{n\beta_\tau}. \quad (3.16)$$

This proves (3.6). □

***Remark 3.6.*** In the proof of Lemma 3.5, if we apply the generalized Hölder inequality to the first integral in (3.8) with $N$ positive numbers $p_1, \ldots, p_N$ defined by

$$p_\ell = \sum_{i=1}^N \frac{\overline{H}_\ell}{\overline{H}_i} \qquad (\ell = 1, \ldots, N),$$

then the above proof shows that if $T \subseteq I \in \mathcal{A}$, then, similarly to (3.15), the following inequality holds:

$$\mathbb{E}[L(x,T)^n] \leq c_{3,10}^n (n!)^{N\nu} \lambda_N(T)^{n(1-\nu)}, \quad (3.17)$$

where $\nu = d/(\sum_{\ell=1}^N \frac{1}{\overline{H}_\ell}) \in (0,1)$. We will apply this inequality in the proof of Theorem 3.4 below.

**Lemma 3.7.** *Assume the conditions of Theorem 3.4 hold. Let $\tau$ be the unique integer in $\{1,\ldots,N\}$ satisfying (3.5). There then exists a positive and finite constant $c_n$, depending on $N, d, H, I$ and $n$, such that for all subintervals $T = [a, a + \langle h \rangle] \subseteq I$ with $h > 0$ small, $x, y \in \mathbb{R}^d$ with $|x - y| \leq 1$, all even integers $n \geq 1$ and all $\gamma \in (0, 1 \wedge \frac{\alpha_\tau}{2\tau})$,*

$$\mathbb{E}[(L(x,T) - L(y,T))^n] \leq c_n |x-y|^{n\gamma} h^{n(\beta_\tau - \overline{H}_\tau \gamma)}, \tag{3.18}$$

*where $\alpha_\tau = \sum_{\ell=1}^\tau \frac{1}{\overline{H}_\ell} - d$.*

**Proof.** Let $\gamma \in (0, 1 \wedge \frac{\alpha_\tau}{2\tau})$ be a constant. Note that by the elementary inequalities

$$|e^{iu} - 1| \leq 2^{1-\gamma} |u|^\gamma \quad \text{for all } u \in \mathbb{R} \tag{3.19}$$

and $|u + v|^\gamma \leq |u|^\gamma + |v|^\gamma$, we see that for all $u^1, \ldots, u^n, x, y \in \mathbb{R}^d$,

$$\prod_{j=1}^n |e^{-i\langle u^j, x\rangle} - e^{-i\langle u^j, y\rangle}| \leq 2^{(1-\gamma)n} |x-y|^{n\gamma} {\sum}' \prod_{j=1}^n |u_{k_j}^j|^\gamma, \tag{3.20}$$

where the summation $\sum'$ is taken over all the sequences $(k_1, \ldots, k_n) \in \{1, \ldots, d\}^n$.

It follows from (2.44), (3.20) and Lemma 2.6 that for every even integer $n \geq 2$,

$$\begin{aligned}
&\mathbb{E}[(L(x,T) - L(y,T))^n] \\
&\leq |x-y|^{n\gamma} {\sum}' \int_{T^n} d\overline{t} \int_{\mathbb{R}^{nd}} \prod_{m=1}^n |u_{k_m}^m|^\gamma \exp\left[-\tfrac{1}{2} \operatorname{Var}\left(\sum_{j=1}^n \langle u^j, B^{H(t^j)}(t^j)\rangle\right)\right] d\overline{u} \\
&= |x-y|^{n\gamma} {\sum}' \int_{T^n} d\overline{t} \int_{\mathbb{R}^{nd}} \prod_{m=1}^n |u_{k_m}^m|^\gamma \prod_{k=1}^d \exp\left[-\tfrac{1}{2} \operatorname{Var}\left(\sum_{j=1}^n u_k^j B_k^{H(t^j)}(t^j)\right)\right] d\overline{u} \\
&\leq |x-y|^{n\gamma} {\sum}' \int_{T^n} d\overline{t} \int_{\mathbb{R}^{nd}} \prod_{m=1}^n |u_{k_m}^m|^\gamma \prod_{k=1}^d \exp\left[-\tfrac{1}{2} \sum_{\ell=1}^N \operatorname{Var}\left(\sum_{j=1}^n u_k^j Y_\ell(t^j)\right)\right] d\overline{u} \\
&\leq |x-y|^{n\gamma} {\sum}' \int_{T^n} d\overline{t} \int_{\mathbb{R}^{nd}} \prod_{m=1}^n |u_{k_m}^m|^\gamma \prod_{k=1}^d \exp\left[-\tfrac{1}{2} \sum_{\ell=1}^\tau \operatorname{Var}\left(\sum_{j=1}^n u_k^j Y_\ell(t^j)\right)\right] d\overline{u} \\
&= |x-y|^{n\gamma} {\sum}' \int_{T^n} d\overline{t} \prod_{k=1}^d \int_{\mathbb{R}^n} \prod_{j=1}^n |u_k^j|^{\gamma \eta_k^j} \exp\left[-\tfrac{1}{2} \sum_{\ell=1}^\tau \operatorname{Var}\left(\sum_{j=1}^n u_k^j Y_\ell(t^j)\right)\right] dU_k,
\end{aligned} \tag{3.21}$$

where $\eta_k^j = 1$ if $k = k_j$ and $\eta_k^j = 0$ otherwise. Note that for every $j \in \{1, \ldots, n\}$, we have $\sum_{k=1}^d \eta_k^j = 1$.

We take $\beta_\ell = \overline{H}_\ell$ $(1 \leq \ell \leq N)$, $\Delta = 1/n$ and $q = d$ in Lemma 2.10 and let $p_\ell$ $(\ell = 1, \ldots, \tau)$ be the constants satisfying (2.36) and (2.37). Observe that since $\gamma \in (0, \frac{\alpha_\tau}{2\tau})$, there exists an $\ell_0 \in \{1, \ldots, \tau\}$ such that

$$\frac{\overline{H}_{\ell_0} d}{p_{\ell_0}} + 2\overline{H}_{\ell_0}\gamma < 1. \tag{3.22}$$

Combining (3.21) with the generalized Hölder inequality, we have that

$$\begin{aligned}
&\mathbb{E}[(L(x,T) - L(y,T))^n] \\
&\leq |x-y|^{n\gamma} {\sum}' \int_{T^n} \mathrm{d}\bar{t} \\
&\quad \times \prod_{k=1}^{d} \left\{ \left( \int_{\mathbb{R}^n} \prod_{j=1}^{n} |u_k^j|^{\gamma \eta_k^j p_{\ell_0}} \exp\left[-\tfrac{1}{2}\mathrm{Var}\left(\sum_{j=1}^{n} u_k^j Y_{\ell_0}(t^j)\right)\right] \mathrm{d}U_k \right)^{1/p_{\ell_0}} \right. \\
&\quad \left. \times \prod_{\ell \neq \ell_0}^{\tau} \left( \int_{\mathbb{R}^n} \exp\left[-\tfrac{1}{2}\mathrm{Var}\left(\sum_{j=1}^{n} u_k^j Y_{\ell}(t^j)\right)\right] \mathrm{d}U_k \right)^{1/p_\ell} \right\}.
\end{aligned} \tag{3.23}$$

For any $n$ points $t^1, \ldots, t^n \in T$, let $\pi_1, \ldots, \pi_N$ be $N$ permutations of $\{1, 2, \ldots, n\}$ such that for every $1 \leq \ell \leq N$,

$$t_\ell^{\pi_\ell(1)} \leq t_\ell^{\pi_\ell(2)} \leq \cdots \leq t_\ell^{\pi_\ell(n)}. \tag{3.24}$$

Let

$$\mathcal{M}_{\ell_0} = \int_{\mathbb{R}^n} \prod_{j=1}^{n} |u_k^j|^{\gamma \eta_k^j p_{\ell_0}} \exp\left[-\tfrac{1}{2}\mathrm{Var}\left(\sum_{j=1}^{n} u_k^j Y_{\ell_0}(t^j)\right)\right] \mathrm{d}U_k. \tag{3.25}$$

By changing the variables of the above integral by means of the transformation

$$u_k^{\pi_{\ell_0}(j)} = v_k^j - v_k^{j+1}, \qquad j = 1, \ldots, n; \qquad u_k^{\pi_{\ell_0}(n)} = v_k^n,$$

we have that

$$\sum_{j=1}^{n} u_k^j Y_{\ell_0}(t^j) = \sum_{j=1}^{n} v_k^j (Y_{\ell_0}(t^{\pi_{\ell_0}(j)}) - Y_{\ell_0}(t^{\pi_{\ell_0}(j-1)})),$$

where $t^{\pi_{\ell_0}(0)} = 0$.

Furthermore, by the elementary inequality that for $\xi > 0$, $|a-b|^\xi \leq c_\xi(|a|^\xi + |b|^\xi)$, where $c_\xi = 2^{\xi-1}$ if $\xi > 1$ and $1$ if $\xi \leq 1$, we have that

$$\begin{aligned}
\prod_{j=1}^{n} |u_k^{\pi_{\ell_0}(j)}|^{\gamma \eta_k^j p_{\ell_0}} &= \prod_{j=1}^{n-1} |v_k^j - v_k^{j+1}|^{\gamma \eta_k^j p_{\ell_0}} |v_k^n|^{\gamma \eta_k^n p_{\ell_0}} \\
&\leq c^n \prod_{j=1}^{n-1} (|v_k^j|^{\gamma \eta_k^j p_{\ell_0}} + |v_k^{j+1}|^{\gamma \eta_k^j p_{\ell_0}}) |v_k^n|^{\gamma \eta_k^n p_{\ell_0}}.
\end{aligned} \tag{3.26}$$

Moreover, the last product is equal to a finite sum of terms, each of the form $\prod_{j=1}^{n} |v_k^j|^{\gamma \eta_k^j p_{\ell_0} \varepsilon_j}$, where $\varepsilon_j = 0, 1,$ or 2 and $\sum_{j=1}^{n} \sum_{k=1}^{d} \eta_k^j \varepsilon_j = n$.

Let $\sigma_{\ell_0,j}^2 = \mathbb{E}[(Y_{\ell_0}(t^{\pi_{\ell_0}(j)}) - Y_{\ell_0}(t^{\pi_{\ell_0}(j-1)}))^2]$. By Proposition 2.5, we know that $\mathcal{M}_{\ell_0}$ is dominated by the sum over all possible choices of $(\varepsilon_1, \ldots, \varepsilon_n) \in \{0,1,2\}^n$ of the following terms

$$\int_{\mathbb{R}^n} \prod_{j=1}^{n} |v_k^j|^{\gamma \eta_k^j p_{\ell_0} \varepsilon_j} \exp\left(-\frac{C_n}{2} \sum_{j=1}^{n} (v_k^j)^2 \sigma_{\ell_0,j}^2\right) dV_k, \tag{3.27}$$

where $V_k = (v_k^1, \ldots, v_k^n) \in \mathbb{R}^n$. By another change of variable, $w_k^j = \sigma_{\ell_0,j} v_k^j$, the integral in (3.27) can be represented by

$$\prod_{j=1}^{n} \sigma_{\ell_0,j}^{-1-\gamma \eta_k^j p_{\ell_0} \varepsilon_j} \int_{\mathbb{R}^n} \prod_{j=1}^{n} |w_k^j|^{\gamma \eta_k^j p_{\ell_0} \varepsilon_j} \exp\left(-\frac{C_n}{2} \sum_{j=1}^{n} (w_k^j)^2\right) dW_k$$
$$:= C_{n,1} \prod_{j=1}^{n} \sigma_{\ell_0,j}^{-1-\gamma \eta_k^j p_{\ell_0} \varepsilon_j}, \tag{3.28}$$

where

$$C_{n,1} = \int_{\mathbb{R}^n} \prod_{j=1}^{n} |w_k^j|^{\gamma \eta_k^j p_{\ell_0} \varepsilon_j} \exp\left(-\frac{C_n}{2} \sum_{j=1}^{n} (w_k^j)^2\right) dW_k$$

is a constant depending on $n$. Thus, we have obtained that

$$\mathcal{M}_{\ell_0} \leq c_n \prod_{j=1}^{n} \sigma_{\ell_0,j}^{-1-\gamma \eta_k^j p_{\ell_0} \varepsilon_j}. \tag{3.29}$$

The other integrals for $\ell \neq \ell_0$ in (3.23) are easier and can be estimated similarly.

Combining (3.23) with Lemma 2.4 (which gives lower bounds for $\sigma_{\ell,j}^2$), (3.28), (3.29) and the definition of $T$, we have

$$\mathbb{E}[(L(x,T) - L(y,T))^n]$$
$$\leq c_n |x-y|^{n\gamma} {\sum}' \left\{ \int_{\Pi_{\ell_0}} \prod_{j=1}^{n} (t_{\ell_0}^{\pi_{\ell_0}(j)} - t_{\ell_0}^{\pi_{\ell_0}(j-1)})^{-(\overline{H}_{\ell_0} d/p_{\ell_0}) - \gamma \overline{H}_{\ell_0} \varepsilon_j} dt_{\ell_0}^1 \cdots dt_{\ell_0}^n \right. \tag{3.30}$$
$$\left. \times \prod_{\ell \neq \ell_0}^{\tau} \int_{\Pi_\ell} \prod_{j=1}^{n} \frac{1}{(t_\ell^{\pi_\ell(j)} - t_\ell^{\pi_\ell(j-1)})^{\overline{H}_\ell d/p_\ell}} dt_\ell^1 \cdots dt_\ell^n \times \prod_{\ell=\tau+1}^{N} h_\ell^n \right\}.$$

In the above, $\Pi_\ell = \{a_\ell \leq t_\ell^{\pi_\ell(1)} \leq \cdots \leq t_\ell^{\pi_\ell(n)} \leq a_\ell + h_\ell\}$ for every $1 \leq \ell \leq \tau$.

We now take $h_1 = \cdots = h_N = h < \delta$ in (3.30). Then by Lemma 2.11 and noting that (3.22) holds, we have

$$\mathbb{E}[(L(x,T) - L(y,T))^n] \leq c_n |x-y|^{n\gamma} h^{n(N-(1-\frac{1}{n}))\sum_{\ell=1}^{\tau} \overline{H}_\ell d/p_\ell - \overline{H}_{\ell_0}\gamma)}$$
$$\leq c_n |x-y|^{n\gamma} h^{n(\beta_\tau - \overline{H}_\tau \gamma)}, \quad (3.31)$$

where the last inequality follows from Lemma 2.10 and the fact that $\overline{H}_{\ell_0} \leq \overline{H}_\tau$. □

The proof of Theorem 3.4 is similar to the proofs of Xiao and Zhang [45] Theorem 4.1 and Ayache, Wu and Xiao [5] Theorem 3.1; we include it here for completeness.

**Proof Theorem 3.4.** Let $I = [a,b] \in \mathcal{A}$ be fixed and assume that (3.5) holds. It follows from Lemma 3.7 and the multiparameter version of Kolmogorov's continuity theorem (cf. Khoshnevisan [25]) that for every $T \subseteq I$, the mfBs $\{B^{H(t)}\}$ has almost surely a local time $L(x,T)$ that is continuous for all $x \in \mathbb{R}^d$.

To prove the joint continuity, observe that for all $x, y \in \mathbb{R}^d$ and $s, t \in I$ with $|s-t| < \delta$, where $\delta > 0$ is the same as in Lemma 2.2, we have

$$\mathbb{E}[(L(x,[a,s]) - L(y,[a,t]))^n]$$
$$\leq 2^{n-1}\{\mathbb{E}[(L(x,[a,s]) - L(x,[a,t]))^n] + \mathbb{E}[(L(x,[a,t]) - L(y,[a,t]))^n]\}. \quad (3.32)$$

Since the difference $L(x,[a,s]) - L(x,[a,t])$ can be written as a sum of a finite number (which only depends on $N$) of terms of the form $L(x,T_j)$, where each $T_j \in \mathcal{A}$ is a closed subinterval of $I$ with at least one edge length $\leq |s-t|$, we can use Lemma 3.5 and Remark 3.6 to bound the first term in (3.32). On the other hand, the second term in (3.32) can be dealt with using Lemma 3.7, as above. Consequently, for some small $\gamma \in (0,1)$, the right-hand side of (3.32) is bounded by $c_{3,11}^n (|x-y| + |s-t|)^{n\gamma}$, where $n \geq 2$ is an arbitrary even integer. Therefore, the joint continuity of the local times on $I$ again follows from the multiparameter version of Kolmogorov's continuity theorem. This finishes the proof of Theorem 3.4. □

Similarly to the proof of Corollary 3.2, we derive from Theorem 3.4 the following, more general, result.

**Corollary 3.8.** *Let $I \in \mathcal{A}$ be fixed. If $\sum_{\ell=1}^{N} \frac{1}{H_\ell(t)} > d$ for all $t \in I$, then mfBs $\{B^{H(t)}(t)\}$ has a jointly continuous local time on $I$ almost surely.*

The next corollary is a direct consequence of Corollary 3.8 and the continuity of the Hurst functionals. We state it to emphasize that the joint continuity of the local time is a local property depending on $H(t)$.

**Corollary 3.9.** *Let $t^0 \in [\varepsilon, \Lambda]^N$ be fixed. If $\sum_{\ell=1}^{N} \frac{1}{H_\ell(t^0)} > d$, then mfBs $\{B^{H(t)}(t)\}$ has a jointly continuous local time on $U(t^0, r_0)$ for some $r_0 > 0$ almost surely, where $U(t^0, r_0)$ is the open ball centered at $t^0$ with radius $r_0$.*

**Remark 3.10.** Let $\{\widehat{B}_0^{H(t)}(t), t \in \mathbb{R}_+^N\}$ be a real-valued multifractional Brownian sheet defined by the harmonizable representation (cf. Herbin [22])

$$\widehat{B}_0^{H(t)}(t) = \int_{\mathbb{R}^N} \prod_{\ell=1}^N \frac{e^{it_\ell \lambda_\ell} - 1}{|\lambda_\ell|^{H_\ell(t)+1/2}} \widehat{W}(d\lambda), \forall t \in \mathbb{R}_+^N, \quad (3.33)$$

where $\widehat{W}$ is the complex-valued Gaussian random measure in $\mathbb{R}^N$ with the Lebesgue measure $\Lambda_N$ as its control measure. Based on Dobríc and Ojeda [16] (see also Stoev and Taqqu [35]), $\{\widehat{B}_0^{H(t)}(t)\}$ and $\{{B'}_0^{H(t)}(t)\}$ defined by equation (1.4) are equivalent up to a multiplicative deterministic function. For $\{\widehat{B}_0^{H(t)}(t)\}$, we can prove a similar result as in Lemma 2.2 (cf. Boufoussi, Dozzi and Guerbaz [13], Wu [38]). Hence, the existence of local times of the $(N,d)$-multifractional Brownian sheet $\{\widehat{B}^{H(t)}(t)\}$ defined by (1.5), where $\{B_1^{H(t)}(t)\}, \ldots, \{B_d^{H(t)}(t)\}$ are independent copies of $\{\widehat{B}_0^{H(t)}(t)\}$, can be proven as in Theorem 3.1 and Corollary 3.2. However, we have not been able to prove that $\{\widehat{B}_0^{H(t)}(t)\}$ is (even one-sided) sectorial locally non-deterministic (since the Fourier analytic technique employed in Wu and Xiao [39] for proving that fractional Brownian sheets are sectorial local non-deterministic fails when $H(t)$ varies in $t$). It is an open problem to establish the joint continuity of the local times of $\{\widehat{B}^{H(t)}(t)\}$.

## 4. Hölder conditions for $L(x, \bullet)$

In this section, we investigate the local and global asymptotic behavior of the local time $L(x, \cdot)$ at $x$ as a random measure on $\mathbb{R}_+^N$. Results in this section carry information about fractal properties of the sample functions of mfBs; see Section 5.

By applying Lemma 3.5, we can prove the following technical lemmas, which will be useful in this section.

**Lemma 4.1.** *Under the conditions of Lemma 3.5, there exists a positive and finite constant $c_{4,1}$, depending only on $N, d$, $H$ and $I$, such that for all $a \in I$ and hypercubes $T = [a, a + \langle r \rangle] \subseteq I$ with $r < \delta$, where $\delta > 0$ is the same as in Lemma 2.2, $x \in \mathbb{R}^d$ and all integers $n \geq 1$,*

$$\mathbb{E}[L(x + B^{H(a)}(a), T)^n] \leq c_{4,1}^n (n!)^{N-\beta_\tau} r^{n\beta_\tau}. \quad (4.1)$$

**Proof.** The proof is similar to the proof of Ayache, Wu and Xiao [5] Lemma 3.11 and we include it here for completeness. For each fixed $a \in I$, we define the Gaussian random field $Y = \{Y(t), t \in \mathbb{R}_+^N\}$ with values in $\mathbb{R}^d$ by $Y(t) = B^{H(t)}(t) - B^{H(a)}(a)$. It follows from (2.41) that if $\{B^{H(t)}(t)\}$ has a local time $L(x, S)$ on any Borel set $S$, then $Y$ also has a local time $\tilde{L}(x, S)$ on $S$ and, moreover, $L(x + B^{H(a)}(a), S) = \tilde{L}(x, S)$. With a little modification, the proof of Lemma 3.5 works for the Gaussian field $Y$. Hence, we derive that (4.1) holds. □

The following lemma is a consequence of Lemma 4.1 and Chebyshev's inequality.

**Lemma 4.2.** *Under the conditions of Lemma 3.5, there exist positive constants $c_{4,2}$, $b$ (depending only on $N$, $d$, $I$ and $H$) such that for all $a \in I$, $T = [a, a + \langle r \rangle]$ with $r \in (0, \delta)$, where $\delta > 0$ is the same as in Lemma 2.2, $x \in \mathbb{R}^d$ and $u > 1$ large enough, we have*

$$\mathbb{P}\{L(x + B^{H(a)}(a), T) \geq c_{4,2} r^{\beta_\tau} u^{N-\beta_\tau}\} \leq \exp(-bu). \tag{4.2}$$

Let $U(t, r)$ be the open ball centered at $t$ with radius $r$, let $H^U$ be the vector defined by (3.1) with $I = U(t, r)$ and let $\tau(U)$ be the positive integer satisfying the corresponding condition (3.5). By applying Lemma 3.5 and the Borel–Cantelli lemma, one can easily derive the following law of the iterated logarithm for the local time $L(x, \cdot)$: there exists a positive constant $c_{4,3}$ such that for every $x \in \mathbb{R}^d$ and $t \in (0, \infty)^N$,

$$\limsup_{r \to 0} \frac{L(x, U(t, r))}{\varphi_U(r)} \leq c_{4,3}, \tag{4.3}$$

where $\varphi_U(r) = r^{\beta_{\tau(U)}} (\log \log(1/r))^{N - \beta_{\tau(U)}}$ with

$$\beta_{\tau(U)} = N - \tau(U) - H^U_{\tau(U)} d + \sum_{\ell=1}^{\tau(U)} \frac{H^U_{\tau(U)}}{H^U_\ell}.$$

Because of the continuity of $H_\ell(t)$ $(1 \leq \ell \leq N)$, it can be verified that

$$\tau(U) \to \tau(t) \quad \text{and} \quad \beta_{\tau(U)} \to \beta_{\tau(t)} \quad \text{as } r \to 0, \tag{4.4}$$

where $\tau(t)$ is the unique integer satisfying

$$\sum_{\ell=1}^{\tau(t)-1} \frac{1}{H_\ell(t)} \leq d < \sum_{\ell=1}^{\tau(t)} \frac{1}{H_\ell(t)} \tag{4.5}$$

and

$$\beta_{\tau(t)} = N - \tau(t) - H_{\tau(t)}(t) d + \sum_{\ell=1}^{\tau(t)} \frac{H_{\tau(t)}(t)}{H_\ell(t)}. \tag{4.6}$$

It follows from Fubini's theorem that with probability one, (4.3) holds for almost all $t \in (0, \infty)^N$. We now prove a stronger version of this result, which is useful in determining the Hausdorff dimension of the level set.

**Theorem 4.3.** *Let $I \in \mathcal{A}$ be a fixed interval and assume that $d < \sum_{\ell=1}^N \frac{1}{H_\ell}$. For any fixed $x \in \mathbb{R}^d$, let $L(x, \cdot)$ be the local time of $\{B^{H(t)}(t)\}$ at $x$ which is a random measure supported on the level set. There then exists a positive and finite constant $c_{4,4}$, independent of $x$, such that with probability 1, the following holds for $L(x, \cdot)$-almost all $t \in I$:*

$$\limsup_{r \to 0} \frac{L(x, U(t, r))}{\varphi_t(r)} \leq c_{4,4}, \tag{4.7}$$

where $\varphi_t(r) = r^{\beta_{\tau(t)}}(\log\log(1/r))^{N-\beta_{\tau(t)}}$ and where $\beta_{\tau(t)}$ is defined by (4.6).

**Proof.** For every integer $k > 0$, we consider the random measure $L_k(x, \bullet)$ on the Borel subsets $C$ of $I$ defined by

$$L_k(x, C) = \int_C (2\pi k)^{d/2} \exp\left(-\frac{k|B^{H(t)}(t) - x|^2}{2}\right) dt$$
$$= \int_C \int_{\mathbb{R}^d} \exp\left(-\frac{|\xi|^2}{2k} + i\langle \xi, B^{H(t)}(t) - x\rangle\right) d\xi\, dt. \quad (4.8)$$

Then, by the occupation density formula (2.41) and the continuity of the function $y \mapsto L(y, C)$, one can verify that almost surely $L_k(x, C) \to L(x, C)$ as $k \to \infty$ for every Borel set $C \subseteq I$.

For every integer $m \geq 1$, let $f_m(t) = L(x, U(t, 2^{-m}))$. From the proof of Theorem 3.4, we can see that almost surely the functions $f_m(t)$ are continuous and bounded. Hence, we have, almost surely, for all integers $m, n \geq 1$,

$$\int_I [f_m(t)]^n L(x, dt) = \lim_{k \to \infty} \int_I [f_m(t)]^n L_k(x, dt). \quad (4.9)$$

It follows from (4.9), (4.8) and the proof of Pitt [34], Proposition 3.1 that for every positive integer $n \geq 1$,

$$\mathbb{E} \int_I [f_m(t)]^n L(x, dt)$$
$$= \left(\frac{1}{2\pi}\right)^{(n+1)d} \int_I \int_{U(t, 2^{-m})^n} \int_{\mathbb{R}^{(n+1)d}} \exp\left(-i \sum_{j=1}^{n+1} \langle x, u^j\rangle\right) \quad (4.10)$$
$$\times \mathbb{E} \exp\left(i \sum_{j=1}^{n+1} \langle u^j, B^{H(s^j)}(s^j)\rangle\right) d\overline{u}\, d\overline{s},$$

where $\overline{u} = (u^1, \ldots, u^{n+1}) \in \mathbb{R}^{(n+1)d}$ and $\overline{s} = (t, s^1, \ldots, s^n)$. Similarly to the proof of (3.6), for sufficiently large $m$, we have that the right-hand side of equation (4.10) is at most

$$c_{4,3}^n \int_I \int_{U(t, 2^{-m})^n} \frac{d\overline{s}}{[\det \text{Cov}(B_0^{H(t)}(t), B_0^{H(s^1)}(s^1), \ldots, B_0^{H(s^n)}(s^n))]^{d/2}}$$
$$\leq c_{4,4}^n (n!)^{N-\beta_\tau(U)} 2^{-mn\beta_\tau(U)}, \quad (4.11)$$

where $c_{4,4}$ is a positive finite constant depending only on $N, d, H$, and $I$.

Let $\gamma > 0$ be a constant, the value of which will be determined later. We consider the random set

$$I_m(\omega) = \{t \in I : f_m(t) \geq \gamma \varphi_U(2^{-m})\}.$$

Denote by $\mu_\omega$ the restriction of the random measure $L(x, \cdot)$ to $I$, that is, $\mu_\omega(E) = L(x, E \cap I)$ for all Borel sets $E \in \mathbb{R}_+^N$. We now take $n = \lfloor \log m \rfloor$, where $\lfloor y \rfloor$ denotes the integer part of $y$. Then by applying (4.11) and by Stirling's formula, we have

$$\mathbb{E}\mu_\omega(I_m) \leq \frac{\mathbb{E}\int_I [f_m(t)]^n L(x, \mathrm{d}t)}{[\gamma \varphi_U(2^{-m})]^n} \qquad (4.12)$$
$$\leq \frac{c_{4,4}^n (n!)^{N-\beta_\tau(U)} 2^{-mn\beta_\tau(U)}}{\gamma^n 2^{-mn\beta_\tau(U)} (\log m)^{n(N-\beta_\tau(U))}} \leq m^{-2},$$

provided $\gamma > 0$ is chosen large enough, say, $\gamma \geq c_{4,2}$. This implies that

$$\mathbb{E}\left(\sum_{m=1}^\infty \mu_\omega(I_m)\right) < \infty.$$

Therefore, with probability 1 for $\mu_\omega$-almost all $t \in I$, we have

$$\limsup_{m \to \infty} \frac{L(x, U(t, 2^{-m}))}{\varphi_U(2^{-m})} \leq c_{4,2}. \qquad (4.13)$$

By (4.4), we can see that for $m$ sufficiently large, there exists a constant $c_{4,5} > 0$ such that $\varphi_U(2^{-m}) \leq c_{4,5} \varphi_t(2^{-m})$. Therefore, we have

$$\limsup_{m \to \infty} \frac{L(x, U(t, 2^{-m}))}{\varphi_t(2^{-m})} \leq c_{4,6}. \qquad (4.14)$$

Finally, for any $r > 0$ small enough, there exists an integer $m$ such that $2^{-m} \leq r < 2^{-m+1}$ and (4.13) is applicable. Since $\varphi_t$ is increasing near 0, (4.7) follows from (4.13) and a monotonicity argument. □

Recall that the pointwise Hölder exponent of a random field $\{X(t), t \in \mathbb{R}^N\}$ at a point $t \in \mathbb{R}^N$ is defined by

$$\alpha_X(t) = \sup\left\{\alpha : \lim_{|h| \to 0} \frac{X(t+h) - X(t)}{|h|^\alpha} = 0\right\}. \qquad (4.15)$$

By Theorem 4.3, and noting that $L(\cdot, t)$ vanishes outside some compact set $U$ (depending on $\omega$), we have the following corollary.

**Corollary 4.4.** *For every $x \in \mathbb{R}^d$, the pointwise Hölder exponent $\alpha_L$ of $L(x, t)$ at $t$ satisfies*

$$\alpha_L(t) \geq \beta_{\tau(t)} \qquad a.s. \qquad (4.16)$$

## 5. Hausdorff dimensions of the level sets

For $x \in \mathbb{R}^d$, let $\Gamma_x = \{t \in (0,\infty)^N : B^{H(t)}(t) = x\}$ be the level set of the multifractional Brownian sheet $\{B^{H(t)}(t)\}$. In this section, we determine the Hausdorff dimension of $\Gamma_x$. We remark that the corresponding problem of finding the Hausdorff dimensions of the level sets of multifractional Brownian motion has been investigated by Boufoussi *et al.* [13] and the Hausdorff dimensions of the level sets of fractional Brownian sheets were studied in Ayache and Xiao [6]. As shown by the following theorem, the fractal structure of $\Gamma_x$ is much richer than the level sets of multifractional Brownian motion and is locally reminiscent of the level sets of a fractional Brownian sheet.

**Theorem 5.1.** *Let $\{B^{H(t)}(t)\} = \{B^{H(t)}(t), t \in \mathbb{R}_+^N\}$ be an $(N,d)$-multifractional Brownian sheet with Hurst functionals $H_\ell(t)$ $(\ell = 1, \ldots, N)$. For any interval $I \in \mathcal{A}$, let $t^* \in I$ be a point satisfying*

$$\sum_{\ell=1}^{N} \frac{1}{H_\ell(t^*)} = \max_{t \in I} \left\{ \sum_{\ell=1}^{N} \frac{1}{H_\ell(t)} \right\}$$

*and*

$$0 < H_1(t^*) \leq \cdots \leq H_N(t^*) < 1.$$

*If $\sum_{\ell=1}^{N} \frac{1}{H_\ell(t^*)} < d$, then for every $x \in \mathbb{R}^d$, we have $\Gamma_x \cap I = \varnothing$ a.s. If $\sum_{\ell=1}^{N} \frac{1}{H_\ell(t^*)} > d$, then for any $x \in \mathbb{R}^d$, with positive probability,*

$$\dim_{\mathrm{H}}(\Gamma_x \cap I) = \beta_{\tau(t^*)}, \tag{5.1}$$

*where*

$$\beta_{\tau(t^*)} = \min\left\{ \sum_{\ell=1}^{k} \frac{H_k(t^*)}{H_\ell(t^*)} + N - k - H_k(t^*)d, 1 \leq k \leq N \right\}.$$

**Remark 5.2.** It can be verified that if (4.5) holds for $t = t^*$, then $\beta_{\tau(t^*)}$ is the same as in (4.6) with $t$ replaced by $t^*$. In the special case where $H(t) = H$ is a constant, $\beta_{\tau(t^*)}$ reduces to the form derived in Ayache and Xiao [6] for the Hausdorff dimension of the level sets of a fractional Brownian sheet.

When $\sum_{\ell=1}^{N} \frac{1}{H_\ell(t^*)} = d$, we believe that for every $x \in \mathbb{R}^d$, $\Gamma_x \cap I = \varnothing$ a.s. However, the method of this paper is not sufficient for proving this statement.

**Proof of Theorem 5.1.** We will follow the proof of Ayache and Xiao [6] Theorem 5.

First, we prove

$$\dim_{\mathrm{H}}(\Gamma_x \cap I) \leq \min\left\{ \sum_{\ell=1}^{k} \frac{H_k(t^*)}{H_\ell(t^*)} + N - k - H_k(t^*)d, 1 \leq k \leq N \right\} \quad \text{a.s.} \tag{5.2}$$

and $\Gamma_x \cap I = \varnothing$ a.s. when the right-hand side of (5.2) is negative.

Without loss of any generality, we may assume that $I = [a, a + \langle h \rangle]$ and $h$ is small so that Lemma 2.2 is applicable. For an integer $n \geq 1$, divide the interval $I$ into $m_n \leq n^{\sum_{\ell=1}^{N}(H_\ell(t^*))^{-1}}$ subrectangles $R_{n,j}$ of side lengths $n^{-1/H_\ell(t^*)}h$ ($\ell = 1, \ldots, N$). Let $0 < \varrho < 1$ be fixed and let $\kappa_{n,j}$ be the lower-left vertex of $R_{n,j}$. Define $\rho(s,t) = (\mathbb{E}[B^{H(t)}(t) - B^{H(s)}(s)]^2)^{1/2}$, denote by $N_\rho(R_{n,j}, \varepsilon)$ the smallest number of $\rho$-balls of radius $\varepsilon$ needed to cover $R_{n,j}$ and denote by $D$ the diameter of $R_{n,j}$, that is,

$$D := \sup\{\rho(s,t) : s, t \in R_{n,j}\}.$$

By Lemma 2.2, we have

$$\int_0^D \sqrt{\log N_\rho(R_{n,j}, \varepsilon)}\, d\varepsilon \leq cn^{-1}.$$

Combining the above inequality with Talagrand [36], Lemma 2.1 gives

$$\mathbb{P}\left\{\max_{s,t \in R_{n,j}} |B^{H(s)}(s) - B^{H(t)}(t)| > n^{-(1-\varrho)}\right\} \leq e^{-cn^{2\varrho}}. \tag{5.3}$$

Then for $n$ sufficiently large, we have, for any fixed $x \in I$,

$$\begin{aligned}
&\mathbb{P}\{x \in \{B^{H(u)}(u), u \in R_{n,j}\}\} \\
&\leq \mathbb{P}\left\{\max_{s,t \in R_{n,j}} |B^{H(s)}(s) - B^{H(t)}(t)| \leq n^{-(1-\varrho)}; x \in \{B^{H(u)}(u), u \in R_{n,j}\}\right\} \\
&\quad + \mathbb{P}\left\{\max_{s,t \in R_{n,j}} |B^{H(s)}(s) - B^{H(t)}(t)| > n^{-(1-\varrho)}\right\} \\
&\leq \mathbb{P}\{|B^{H(\kappa_{n,j})}(\kappa_{n,j}) - x| \leq n^{-(1-\varrho)}\} + e^{-cn^{2\varrho}} \\
&\leq c_{5,1} n^{-(1-\varrho)d}.
\end{aligned} \tag{5.4}$$

In the above, we have applied the fact that $\mathrm{Var}(B_0^{H(t)}(t)) \geq c$ for all $t \in I$ to derive the last inequality.

If $\sum_{\ell=1}^{N} \frac{1}{H_\ell(t^*)} < d$, we choose $\varrho > 0$ such that $(1-\varrho)d > \sum_{\ell=1}^{N} \frac{1}{H_\ell(t^*)}$. Let $N_n$ be the number of rectangles $R_{n,j}$ such that $x \in \{B^{H(u)}(u), u \in R_{n,j}\}$. It follows from (5.4) that

$$\mathbb{E}(N_n) \leq c_{5,1} n^{\sum_{\ell=1}^{N}(H_\ell(t^*))^{-1}} n^{-(1-\varrho)d} \to 0 \quad \text{as } n \to \infty. \tag{5.5}$$

Since the random variables $N_n$ are integer-valued, (5.5) and Fatou's lemma imply that a.s. $N_n = 0$ for infinitely many $n$. Therefore, $\Gamma_x \cap I = \varnothing$ a.s.

We now assume $\sum_{\ell=1}^{N} \frac{1}{H_\ell(t^*)} > d$. For every $k \in \{1, 2, \ldots, N\}$, define

$$\eta_k = \sum_{\ell=1}^{k} \frac{H_k(t^*)}{H_\ell(t^*)} + N - k - H_k(t^*)d.$$

By Ayache and Xiao [6], Lemma 7, we have $\eta_k > 0$. Define a covering $\{R'_{n,j}\}$ of $\Gamma_x \cap I$ by $R'_{n,j} = R_{n,j}$ if $x \in \{B^{H(u)}(u), u \in R_{n,j}\}$ and $R'_{n,j} = \varnothing$ otherwise. $R'_{n,j}$ can be covered by $n^{\sum_{\ell=k+1}^{N}(H_k(t^*)^{-1} - H_\ell(t^*)^{-1})}$ cubes of side length $n^{-H_k(t^*)^{-1}}h$. Thus, we can cover the level set $\Gamma_x \cap I$ by a sequence of cubes of side length $n^{-H_k(t^*)^{-1}}h$. Let $\varrho \in (0,1)$ be small and let

$$\eta_k(\varrho) = \sum_{\ell=1}^{k} \frac{H_k(t^*)}{H_\ell(t^*)} + N - k - H_k(t^*)(1-\varrho)d.$$

Clearly, $\eta_k(\varrho) > \eta_k > 0$. In the following, we prove that the Hausdorff dimension of $\Gamma_x \cap I$ is bounded above by $\eta_k(\varrho)$ almost surely. To this end, by (5.4), we have

$$\mathbb{E}\left[\sum_{j=1}^{m_n} n^{\sum_{\ell=k+1}^{N}(H_k(t^*)^{-1} - H_\ell(t^*)^{-1})}(n^{-H_k(t^*)^{-1}})^{\eta_k(\varrho)}\mathbb{1}_{\{x \in \{B^{H(u)}(u), u \in R_{n,j}\}\}}\right]$$
$$\leq c_{5,2} n^{\sum_{\ell=1}^{N}H_\ell(t^*)^{-1} + \sum_{\ell=k+1}^{N}(H_k(t^*)^{-1} - H_\ell(t^*)^{-1}) - \eta_k(\varrho)H_k(t^*)^{-1} - (1-\varrho)d} = c_{5,2}.$$
(5.6)

Fatou's lemma implies that the $\eta_k(\varrho)$-dimensional Hausdorff measure of $\Gamma_x \cap I$ is finite a.s. and thus $\dim_{\mathrm{H}}(\Gamma_x \cap I) \leq \eta_k(\varrho)$ almost surely. Letting $\varrho \downarrow 0$ along the rational numbers, we obtain $\dim_{\mathrm{H}}(\Gamma_x \cap I) \leq \eta_k$ and therefore (5.2).

To prove the lower bound in (5.1), we assume $\tau(t^*) = k$, that is,

$$\sum_{\ell=1}^{k-1} \frac{1}{H_\ell(t^*)} \leq d < \sum_{\ell=1}^{k} \frac{1}{H_\ell(t^*)}.$$

By Condition A, there exists a positive number $\varsigma$ such that for all $t \in I_\varsigma := [t^* - \langle \varsigma \rangle, t^* + \langle \varsigma \rangle] \cap I$, we have

$$\sum_{\ell=1}^{k-1} \frac{1}{\overline{H}_\ell} \leq d < \sum_{\ell=1}^{k} \frac{1}{\overline{H}_\ell}, \tag{5.7}$$

where $\overline{H}_\ell$ ($1 \leq \ell \leq N$) are defined as in (3.1) with $I_\varsigma$ in place of $I$. Note that (5.7) and Ayache and Xiao [6], Lemma 3.3, imply that

$$\sum_{\ell=1}^{k} \frac{\overline{H}_k}{\overline{H}_\ell} + N - k - \overline{H}_k d \in (N-k, N-k+1].$$

Thus, we can choose $\varrho > 0$ such that

$$\gamma := \sum_{\ell=1}^{k} \frac{\overline{H}_k}{\overline{H}_\ell} + N - k - \overline{H}_k(1+\varrho)d > N - k. \tag{5.8}$$

It is sufficient to prove that there is a constant $c_{5,3} > 0$ such that

$$\mathbb{P}\{\dim_{\mathrm{H}}(\Gamma_x \cap I_\varsigma) \geq \gamma\} \geq c_{5,3}. \tag{5.9}$$

Our proof of (5.9) is based on the capacity argument due to Kahane [23]. Similar methods have been used by Adler [1], Testard [37] and Xiao [40], to mention just a few.

Let $\mathcal{M}_\gamma^+$ be the space of all non-negative measures on $\mathbb{R}^N$ with finite $\gamma$-energy. It is known (cf. Adler [1]) that $\mathcal{M}_\gamma^+$ is a complete metric space under the metric

$$\|\mu\|_\gamma = \int_{\mathbb{R}^N} \int_{\mathbb{R}^N} \frac{\mu(\mathrm{d}t)\mu(\mathrm{d}s)}{|t-s|^\gamma}. \tag{5.10}$$

We define a sequence of random positive measures $\mu_n$ on the Borel sets $C \subseteq I_\varsigma$ by

$$\begin{aligned}\mu_n(C) &= \int_C (2\pi n)^{d/2} \exp\left(-\frac{n|B^{H(t)}(t) - x|^2}{2}\right) \mathrm{d}t \\ &= \int_C \int_{\mathbb{R}^d} \exp\left(-\frac{|\xi|^2}{2n} + \mathrm{i}\langle \xi, B^{H(t)}(t) - x\rangle\right) \mathrm{d}\xi\,\mathrm{d}t.\end{aligned} \tag{5.11}$$

It follows from Kahane [23] or Testard [37] that if there are positive and finite constants $c_{5,4}$ and $c_{5,5}$, independent of $\varrho$ and such that

$$\mathbb{E}(\|\mu_n\|) \geq c_{5,4}, \qquad \mathbb{E}(\|\mu_n\|^2) \leq c_{5,5}, \tag{5.12}$$

$$\mathbb{E}(\|\mu_n\|_\gamma) < +\infty, \tag{5.13}$$

where $\|\mu_n\| = \mu_n(I_\varsigma)$, then there is a subsequence of $\{\mu_n\}$, say $\{\mu_{n_k}\}$, such that $\mu_{n_k} \to \mu$ in $\mathcal{M}_\gamma^+$ and $\mu$ is strictly positive with probability $\geq c_{5,4}^2/(2c_{5,5})$. It follows from (4.8) that $\mu$ has its support in $\Gamma_x \cap I_\varsigma$ almost surely. Hence, Frostman's theorem yields (5.9).

It remains to verify (5.12) and (5.13). Let $\sigma^2(t) = \mathrm{Var}(B^{H(t)}(t))$. By Fubini's theorem, we have

$$\begin{aligned}\mathbb{E}(\|\mu_n\|) &= \int_{I_\varsigma} \int_{\mathbb{R}^d} \mathrm{e}^{-\mathrm{i}\langle \xi, x\rangle} \exp\left(-\frac{|\xi|^2}{2n}\right) \mathbb{E}\exp(\mathrm{i}\langle \xi, B^{H(t)}(t)\rangle)\,\mathrm{d}\xi\,\mathrm{d}t \\ &= \int_{I_\varsigma} \int_{\mathbb{R}^d} \mathrm{e}^{-\mathrm{i}\langle \xi, x\rangle} \exp\left(-\frac{1}{2}(n^{-1} + \sigma^2(t))|\xi|^2\right) \mathrm{d}\xi\,\mathrm{d}t \\ &= \int_{I_\varsigma} \left(\frac{2\pi}{n^{-1} + \sigma^2(t)}\right)^{d/2} \exp\left(-\frac{|x|^2}{2(n^{-1} + \sigma^2(t))}\right) \mathrm{d}t \\ &\geq \int_{I_\varsigma} \left(\frac{2\pi}{1 + \sigma^2(t)}\right)^{d/2} \exp\left(-\frac{|x|^2}{2\sigma^2(t)}\right) \mathrm{d}t := c_{5,4}.\end{aligned} \tag{5.14}$$

Denote by $I_{2d}$ the identity matrix of order $2d$, by $\mathrm{Cov}(B^{H(s)}(s), B^{H(t)}(t))$ the covariance matrix of $(B^{H(s)}(s), B^{H(t)}(t))$, let $\Sigma = n^{-1}I_{2d} + \mathrm{Cov}(B^{H(s)}(s), B^{H(t)}(t))$ and denote by

$(\xi, \eta)'$ the transpose of the row vector $(\xi, \eta)$. Then

$$\begin{aligned}
\mathbb{E}(\|\mu_n\|^2) &= \int_{I_\varsigma}\int_{I_\varsigma}\int_{\mathbb{R}^d}\int_{\mathbb{R}^d} e^{-i\langle\xi+\eta,x\rangle} \exp\left(-\frac{1}{2}(\xi,\eta)\Sigma(\xi,\eta)'\right)\,d\xi\,d\eta\,ds\,dt \\
&= \int_{I_\varsigma}\int_{I_\varsigma} \frac{(2\pi)^d}{\sqrt{\det\Sigma}} \exp\left(-\frac{1}{2}(x,x)\Sigma^{-1}(x,x)'\right)\,ds\,dt \qquad (5.15)\\
&\leq \int_{I_\varsigma}\int_{I_\varsigma} \frac{(2\pi)^d}{[\det \mathrm{Cov}(B_0^{H(s)}(s), B_0^{H(t)}(t))]^{d/2}}\,ds\,dt.
\end{aligned}$$

By applying Lemma 2.4 and the regularity of the $H(\cdot)$, it can be proven that for $s, t \in I_\varsigma$,

$$\det \mathrm{Cov}(B_0^{H(s)}(s), B_0^{H(t)}(t)) \geq c_{5,6} \sum_{\ell=1}^{N} |s_\ell - t_\ell|^{2H_\ell(t)} \geq c_{5,6} \sum_{\ell=1}^{N} |s_\ell - t_\ell|^{2\overline{H}_\ell}. \qquad (5.16)$$

Combining (5.7), (5.15), (5.16) and repeatedly applying Lemma 2.8, we obtain

$$\mathbb{E}(\|\mu_n\|^2) \leq c_{5,7} \int_{I_\varsigma}\int_{I_\varsigma} \frac{1}{[\sum_{\ell=1}^{N}|s_\ell-t_\ell|^{2\overline{H}_\ell}]^{d/2}}\,ds\,dt := c_{5,5} < \infty. \qquad (5.17)$$

Thus, we have shown that (5.12) holds.

Similarly to (5.15), we have

$$\begin{aligned}
\mathbb{E}(\|\mu_n\|_\gamma) &= \int_{I_\varsigma}\int_{I_\varsigma} \frac{ds\,dt}{|s-t|^\gamma} \int_{\mathbb{R}^d}\int_{\mathbb{R}^d} e^{-i\langle\xi+\eta,x\rangle} \exp\left(-\frac{1}{2}(\xi,\eta)\Sigma(\xi,\eta)'\right)\,d\xi\,d\eta \\
&\leq c_{5,8} \int_{I_\varsigma}\int_{I_\varsigma} \frac{1}{(\sum_{\ell=1}^{N}|s_\ell-t_\ell|)^\gamma (\sum_{\ell=1}^{N}|s_\ell-t_\ell|^{2\overline{H}_\ell})^{d/2}}\,ds\,dt \qquad (5.18)\\
&\leq c_{5,9} \int_0^\varsigma dt_N \cdots \int_0^\varsigma \frac{1}{(\sum_{\ell=1}^{N} t_\ell^{\overline{H}_\ell})^d (\sum_{\ell=1}^{N} t_\ell)^\gamma}\,dt_1,
\end{aligned}$$

where the two inequalities follow from (5.16) and a change of variables. By using Lemma 2.9 in the same way, we see that $\mathbb{E}(\|\mu_n\|_\gamma) < +\infty$ for any $\gamma$ defined in (5.8). This proves (5.13).

Finally, by letting $\varrho \downarrow 0$, the lower bound for the Hausdorff dimension in (5.1) follows from (5.9) and we have therefore proven Theorem 5.1. $\square$

The proof of Theorem 5.1 suggests that we can consider the Hausdorff dimension of the level set in any neighborhood of a point $t \in (0, \infty)$, provided $\sum_{\ell=1}^{N} \frac{1}{H_\ell(t)} > d$. However, in order to obtain an almost sure result, we have to consider $\Gamma_x$ at a random level $x = B^{H(t)}(t)$. The following corollary can be considered as a local Hausdorff dimension result for the level sets of mfBs.

**Corollary 5.3.** *Let $\{B^{H(t)}(t)\} = \{B^{H(t)}(t), t \in \mathbb{R}_+^N\}$ be an $(N,d)$-multifractional Brownian sheet with Hurst functionals $H_\ell(t)$ ($\ell = 1,\ldots,N$). If $t^0 \in (0,\infty)^N$ satisfies $\sum_{\ell=1}^N \frac{1}{H_\ell(t^0)} > d$, then there exists $r_0 > 0$ such that*

$$\mathbb{P}\left\{\lim_{r \to 0} \dim_{\mathrm{H}}(\Gamma_{B^{H(t)}(t)} \cap U(t^0, r)) = \beta_{\tau(t^0)} \text{ for a.e. } t \in U(t^0, r_0)\right\} = 1, \quad (5.19)$$

*where $\beta_{\tau(t^0)}$ is defined by* (4.6) *with $t^0$ in place of $t$.*

**Proof.** For any $t^0 \in (0,\infty)^N$ such that $\sum_{\ell=1}^N \frac{1}{H_\ell(t^0)} > d$, there exists a positive number $r_0$ such that for all $s \in U(t^0, r_0)$, we have $\sum_{\ell=1}^N \frac{1}{H_\ell(s)} > d$. By Corollary 3.9, mfBs $\{B^{H(t)}(t)\}$ has a jointly continuous local time on $U(t^0, r_0)$.

By (5.2), we have that for every $0 < r < r_0$ and $x \in \mathbb{R}^d$,

$$\dim_{\mathrm{H}}(\Gamma_x \cap U(t^0, r)) \leq \max_{s \in U(t^0, r)} \beta_{\tau(s)} \quad \text{a.s.} \quad (5.20)$$

By (5.20) and Fubini's theorem, we see that

$$\mathbb{P}\left\{\dim_{\mathrm{H}}(\Gamma_x \cap U(t^0, r)) \leq \max_{s \in U(t^0, r)} \beta_{\tau(s)}, \text{ a.e. } x \in \mathbb{R}^d\right\} = 1. \quad (5.21)$$

Since $\{B^{H(t)}(t)\}$ has a local time on $U(t^0, r_0)$, the occupation density formula (2.41) and (5.21) together imply that

$$\mathbb{P}\left\{\dim_{\mathrm{H}}(\Gamma_{B^{H(t)}(t)} \cap U(t^0, r)) \leq \max_{s \in U(t^0, r)} \beta_{\tau(s)} \text{ a.e. } t \in U(t^0, r_0)\right\} = 1. \quad (5.22)$$

On the other hand, by using an argument similar to the one used in the proof of Berman [9] Theorem 2.1 (see also the proof of Xiao [42] Theorem 1.1), we can show that for every $\varepsilon > 0$ and $r \in (0, r_0)$ small enough,

$$\mathbb{P}\left\{\dim_{\mathrm{H}}(\Gamma_{B^{H(t)}(t)} \cap U(t^0, r)) \geq \max_{s \in U(t^0, r)} \beta_{\tau(s)} - \varepsilon \text{ a.e. } t \in U(t^0, r_0)\right\} = 1. \quad (5.23)$$

By letting $r \downarrow 0$ and $\varepsilon \downarrow 0$ along the rational numbers, we see that (5.19) follows from (5.22) and (5.23). □

**Remark 5.4.** Corollary 5.3 shows the explicit way in which the fractal properties of the multifractional Brownian sheet vary in space. In short, the local Hausdorff dimension is derived from the constant parameter formula. Fractional Brownian sheets are essentially fractional integrals of Brownian sheets (cf. Benson *et al.* [8] and Biermé, Meerschaert and Scheffler [12]). The term "multifractional" indicates that the order of fractional integration varies in space. It would be interesting to explore the connection between multifractional Brownian sheets and multifractals. For example, are the level sets of the mfBs multifractals and, if so, how do their structure functions depend on the Hurst index function $H(t)$?

# Acknowledgements

We thank the anonymous referees for their many valuable comments and suggestions. The research of Mark Meerschaert is supported in part by NSF Grant DMS-07-06440. The research of Yimin Xiao is supported in part by NSF Grant DMS-07-06728.